\documentclass[11pt]{article}

\usepackage{amssymb}
\usepackage{amsmath}
\usepackage{theorem}
\usepackage{epsfig}
\usepackage{verbatim}
\usepackage{graphicx}
\usepackage{color}

\newtheorem{thm}{Theorem}[section]

\newtheorem{prop}[thm]{Proposition}

\newtheorem{lemma}[thm]{Lemma}

\newcommand{\R}{\Bbb{R}}

\newcommand{\Z}{\Bbb{Z}}

\newcommand{\T}{\mathbb{T}}
\newcommand{\D}{\displaystyle}

\newcommand{\grad}{\nabla}

\newcommand{\al}{\alpha}
\newcommand{\ep}{\varepsilon}

\newcommand{\g}{\gamma}
\newcommand{\dg}{\partial_{\gamma}}
\newcommand{\dl}{\partial_{\log}}
\newcommand{\p}{\partial}
\newcommand{\e}{\eta}
\newcommand{\de}{\partial_{\eta}}
\newcommand{\dxi}{\partial_{\xi}}
\newcommand{\xt}{\tilde{x}}

\begin{document}

\author{Antonio C\'ordoba, Diego C\'ordoba and Francisco Gancedo}
\title{Uniqueness for SQG patch solutions}

\date{}

\maketitle

\begin{abstract}
 This paper is about the evolution of a temperature front governed by the Surface quasi-geostrophic equation. The existence part of that program within the scale of Sobolev spaces was obtained by one of the authors \cite{FG}. Here we revisit that proof introducing some new tools and points of view which allow us to conclude the also needed uniqueness result.
\end{abstract}

\maketitle

%%%%%%%%%%%%%%%%%%%%%%%%%%%%%%%%%%%%%%%%%%%%%%%%%%%%%%%%%%%%%%%%%%%%%%%%%%%%%%%%%%%%%%%%%%%%%%%%%%%%%%%%%%%%%%%%%%%%%
%%%%%%%%%%%%%%%%%%%%%%%%%%%%%%%%%%%%%%%%%%%%%%%%%%%%%%%%%%%%%%%%%%%%%%%%%%%%%%%%%%%%%%%%%%%%%%%%%%%%%%%%%%%%%%%%%%%%%
%%%%%%%%%%%%%%%%%%%%%%%%%%%%%%%%%%%%%%%%%%%%%%%%%%%%%%%%%%%%%%%%%%%%%%%%%%%%%%%%%%%%%%%%%%%%%%%%%%%%%%%%%%%%%%%%%%%%%

\section{Introduction}

Among the more important partial differential equations of fluid dynamics we have the three dimensional Euler equation, modelling the evolution of an incompressible inviscid fluid, and the surface quasi-geostrophic (SQG) which describes the dynamics of atmospheric temperature \cite{P}. SQG has also the extra mathematical interest of capturing the complexity of 3D Euler equation but in a two dimensional scenario, as was described  in the nowadays classical work \cite{CMT}.   

This model reads
$$
\theta_t+u\cdot\grad \theta=0,
$$
$$
u=(-R_2\theta,R_1\theta),
$$
where $\theta(x,t)$ is the temperature of the 2D fluid with $(x,t)\in \R^2\times[0,+\infty)$. The velocity $u$ is related with the temperature through the Riezs transforms $R_j$ given by
$$
R_j(\theta)(x)=\frac{1}{\pi}\int_{\R^2}\frac{y_j}{|y|^3}\theta(x-y)dy.
$$
Within the equation there is a underlying particle dynamics which preserve the value of $\theta$, implying that the norms $\|\theta\|_{L^p}(t)$, $1\leq p\leq \infty$, remain constants under the evolution.

%Recently there have been an intense work about the equation, see \eqref{}.

In this paper we consider the patch problem, on which the temperature takes two constant values in two complementary domains and the solution of SQG has to be understood in a weak sense, namely: 
\begin{equation}\label{weaksolution}
\int_0^\infty\int_{\R^2}\theta(x,t)(\varphi_t(x,t)+u(x,t)\cdot\nabla\varphi(x,t))dxdt=\int_{\R^2}\theta_0(x)\varphi(x,0)dx,\quad u=(-R_2\theta,R_1\theta),
\end{equation}
for every $\varphi\in C_{c}^\infty([0,\infty)\times\R^2)$. That is, the temperature reads
\begin{equation}\label{patchsolution}
\theta(x,t)=\left\{\begin{array}{rl}
\theta^1,& x\in D^1(t),\\
\theta^2,& x\in D^2(t)=\R^2\setminus D^1(t),
\end{array}\right. 
\end{equation}
where $D^1(t)$ is a simply connected domain.
It gives rise to a contour equation for the free boundary
\begin{equation}\label{parameterization}
\partial D^{j}(t)=\{x(\g,t)=(x_1(\g,t),x_2(\g,t)): \g\in\T\},
\end{equation}
which is moving with the fluid and whose exact formulation can be found in \cite{FG}. It is then clear that the evolution of the patch is equivalent to that of its free boundary $\partial D^{j}(t)$.

This problem was first considered by Resnick in his thesis \cite{Resnick}. Local-in-time existence and uniqueness was proven by Rodrigo \cite{Rodrigo} for $C^\infty$ initial data using Nash-Moser inverse function theorem.  In \cite{FG} the third author proves local-in-time existence for the problem in Sobolev spaces, using energy estimates and properties of a particular parameterization of the contour. Namely, one such that the modulus of the tangent vector to the curve does not depend on the space variable, depending only on time \cite{Hou} and giving us extra cancellations which allows to integrate the system. 
 
In the distributional sense, the gradient of the temperature is given by
$$
\nabla \theta(x,t)=(\theta^2-\theta^1)\dg^{\bot} x(\g,t)\delta(x(\g,t)=x)
$$
for $x(\g,t)$ a given parameterization of the contour and 
$\dg^{\bot} x(\g,t)=(-\dg x_2(\g,t),\dg x_1(\g,t))$. Then Biot-Savart formula helps us to get the velocity field, outside the boundary, in terms of the geometry of the contour, that is
$$
u(x,t)=I_{1}(\nabla^{\bot}\theta)(x,t)=-\frac{\theta^2-\theta^1}{2\pi}\int_{-\pi}^{\pi}\frac{\dg x(\g,t)}{|x-x(\g,t)|}d\g,
$$  
where $I_{1}$ is the Riesz potential of order $1$, which on the Fourier side is  multiplication by $|\xi|^{-1}$. The above integral diverges when $x$ approaches the boundary but only on its tangential component, while its normal component is well defined. This fact is crucial to assign a normal velocity field to the boundary governing its evolution. Since the contribution of the tangential component amount to a reparameterization of the boundary curve, we are free to add such a component satisfying both purposes: to be bounded and having tangent vector with constant length. For a given parameterization $x(\g,t)$, approaching the boundary in both domains we obtain   
$$
u(x(\g,t),t)\cdot \dg^{\bot} x(\g,t)=-\frac{\theta^2-\theta^1}{2\pi}\int_{-\pi}^{\pi}\frac{\dg x(\e,t)\cdot\dg^{\bot} x(\g,t)}{|x(\g,t)-x(\e,t)|}d\e.
$$
And we get the task of finding a good parameterization $x(\g,t)$ and a function $\lambda$ so that  
$$
u(x(\g,t),t)\cdot \dg^{\bot} x(\g,t)=\Big(\frac{\theta^2\!-\!\theta^1}{2\pi}\int_{-\pi}^{\pi}\!\frac{\dg x(\g,t)-\dg x(\e,t)}{|x(\g,t)-x(\e,t)|}d\e+\lambda\dg x(\g,t)\Big)\cdot\dg^{\bot} x(\g,t),
$$ 
and the two purposes mentioned above are achieved. 

Having the length of the vector $\dg x(\g,t)$ as a function in the variable $t$ only  provides the following two identities:
\begin{equation}\label{unacancelacion}
\dg^2 x(\g,t)\cdot\dg x(\g,t)=0,\quad\mbox{and}\quad \dg^3 x(\g,t)\cdot\dg x(\g,t)=-|\dg^2 x(\g,t)|^2.
\end{equation}
The first one gives extra cancellations while the second allows us to perform convenient integration by parts.

Another main character of this play is the so called arc-chord condition which help to control the absence of self-intersections of the boundary curve. This is done through the following quantity:
$$
F(x)(\g,\eta,t)=\frac{|\eta|}{|x(\g,t)-x(\g-\eta,t)|}\quad \forall\,\g,\eta\in[-\pi,\pi],
$$
with $$F(x)(\g,0,t)=\frac{1}{|\partial_\g x(\g,t)|},$$ whose $L^\infty$ norm has to be controlled in the evolution.

Although we can not make justice to the many interesting contribution due to the different authors quoted in our references, let us say that, at the beginning, there was a conjecture about the formation of singularities in the evolution of a vortex patch for Euler equations in dimension two \cite{bertozzi-Majda}. It was disproved by Chemin in a remarkable work \cite{Ch} using paradifferential calculus, and later Bertozzi-Constantin \cite{bertozzi-Constantin} obtained a different proof taking advantage of an extra cancellation satisfied by singular integrals having even kernels.

Between the patch problem for 2D Euler and SQG there is a continuous set of interpolated equations given by
\begin{align}
\begin{split}\label{alpha}
\theta_t+u\cdot\grad \theta=&0,\\
u=(-R_2,R_1)(I_{1-\alpha}\theta)&,\quad 0< \alpha< 1.
\end{split}
\end{align}

The case $\alpha=0$ is the most regular, 2D Euler, while for $\alpha=1$ one gets SQG. The patch problem for those equations was first studied in \cite{CFMR}, where C\'ordoba, Fontelos, Mancho and Rodrigo introduced a very interesting scenario for which they could show numerical evidence of singularity formation: two patches with different temperature approach each other in such a way that they collide at a point where the curvature blows-up. Let us mention that recently it has been shown analytically \cite{GS} that if the curvature is controlled then pointwise collisions can not happen in the patch problem for SQG. In \cite{Scott,ScottDritschel} a different finite time singularity scenario is shown where numerics point at a self-similar blow-up behaviour for SQG patches.  

The system above can also be considered in  more singular cases than SQG, replacing the last identity by the following one
$$
u=(-R_2,R_1)(\Lambda^{\beta}\theta),\quad 0< \beta< 1,
$$
where here $\Lambda=(-\Delta)^{1/2}$, whose Fourier symbol is $|\xi|$. See \cite{CCCGW} for results on this equation with patch solutions.

A classical result in fluid dynamics is the existence for all time of vortex patches for Euler equation which are rotating ellipses \cite{bertozzi-Majda}. The patch problem for the system \eqref{alpha} and SQG present a more complex dynamics, as ellipses are not rotational solutions and some convex interfaces lose this property in finite time \cite{CCGSMZ}. See \cite{Garra} for a study of the growth of the patch support. Recently, in a remarkable series of papers and with an ingenious used of the Crandall-Rabinowitz mountain pass lemma, the authors have extended those global-in-time existence results to a more general class of geometrical shapes for the vortex patch problem \cite{HMV,HMV2}, the $\alpha$-system \eqref{alpha} \cite{HH} and also to the SQG equation \cite{CCG-S,CCG-S2}.  

There are two articles \cite{KYZ,KRYZ} where the patch problem for the $\alpha$-system is considered in a half plane with Dirichlet's condition. The system is proved to be well-posed for $0<\alpha<\frac1{12}$ in the more singular scenario where the patch intersect the fixed boundary. In this framework, singularity formation is shown when two patches of different temperature approach each other.\\

In this paper we will take advantage of a special parameterization of the boundary in the following terms:   

As was mentioned before, patch solutions for the SQG equation are understood in a weak sense. Any such solution with a free boundary given by a smooth parameterization $x(\g,t)$ has to satisfy the equation below
\begin{equation}\label{sqgpatchnormal}
x_t(\g,t)\cdot \dg^{\bot}x(\g,t)=-\int_{\T}\frac{\dg^{\bot}x(\g,t)\cdot\dg x(\g-\e,t)}{|x(\g,t)-x(\g-\e,t)|}d\e,
\end{equation}
where we have taken $\theta_2-\theta_1=\pi$ for the sake of simplicity. On the other hand, any smooth parameterization $x(\g,t)$ satisfying \eqref{sqgpatchnormal} provides a weak SQG solution with the temperature given by (\ref{patchsolution},\ref{parameterization}) (see \cite{FG} for more details).

It is easy to check that the equation above is a reparameterization invariance object, and that the following formula, introduced in \cite{Resnick} and \cite{Rodrigo}, has a well defined tangential velocity and identical normal component
\begin{equation}\label{psqglz}
x_t(\g,t)=\int_{\T}\frac{\dg x(\g,t)-\dg x(\g-\e,t)}{|x(\g,t)-x(\g-\e,t)|}d\e.
\end{equation}
The local-in-time existence result was given in \cite{FG} for initial data satisfying \eqref{unacancelacion} and evolving by  
\begin{align}
\begin{split}\label{QGm}
\D x_t(\g,t)&=\Big(\int_{\T}\frac{\dg x(\g,t)-\dg x(\g-\e,t)}{|x(\g,t)-x(\g-\e,t)|}d\e+\lambda(x)(\g,t)\dg x(\g,t)\Big),\\
\end{split}
\end{align}
\begin{align}
\begin{split}\label{la}
\lambda(x)(\g,t)&=\frac{\g+\pi}{2\pi}\int_\T\frac{\de x(\eta,t)}{|\de x(\eta,t)|^2}\cdot \de \Big(\int_{\T}
\frac{\de x(\eta,t)-\de x(\eta-\xi,t)}{|x(\eta,t)-x(\eta-\xi,t)|}d \xi \Big) d\eta\\
&\quad-\int_{-\pi}^\g \frac{\de x(\e,t)}{|\de x(\e,t)|^2}\cdot \de \Big(\int_{\T}\frac{\de
	x(\e,t)-\de x(\e-\xi,t)}{|x(\e,t)-x(\e-\xi,t)|}d\xi \Big)d\e.
\end{split}
\end{align}
Above $\lambda(x)(-\pi,t)=0$ for the sake of simplicity. In the following we are going to show how it is possible to go from 
(\ref{QGm},\ref{la}) to equation \eqref{psqglz} through a convenient change of variable. This procedure is also valid to go from (\ref{QGm},\ref{la}) to a SQG patch contour equation with a different and more convenient tangential term.  

We denote by $x(\g,t)\in C([0,T];H^3)$ a solution of (\ref{QGm},\ref{la}) and let  $\tilde{x}(\xi,t)$ be given by
\begin{align*}
\tilde{x}(\xi,t)=x(\phi^{-1}(\xi,t),t),\qquad\qquad \g=\phi^{-1}(\xi,t) ,
\end{align*}
or equivalently 
$$
x(\g,t)=\tilde{x}(\phi(\g,t),t),\qquad \xi=\phi(\g,t),
$$
where 
\begin{equation}\label{serunapara}
\phi(\g,t):\R\times\R^+\to\R,\quad \dg\phi(\g,t)>0,\quad \phi(\g,t)-\g\,\,\,\,\mbox{$2\pi$-periodic},
\end{equation}
is a reparameterization in $\g$ for any positive time. Here $\phi$ is a solution of the linear system
\begin{equation}\label{para}
\phi_t(\g,t)=\int_{\T}\frac{\dg\phi(\g,t)-\dg\phi(\e,t)}{|x(\g,t)-x(\e,t)|}d\e+\lambda(x)(\g,t)\dg \phi(\g,t).
\end{equation}
The existence and uniqueness for that system is given in the following proposition, for whose formulation we introduce the space:
$$
H^{\frac{k}{\log}}\equiv\{f\in L^2(\T): \sum_{n\in\Z}\frac{|n|^{2k}}{\log^2(|n|+e)}|\hat{f}(n)|^2=\|f\|^2_{H^{\frac{k}{\log}}}<\infty\}.
$$
\begin{prop}\label{proposicion}
Let $\phi_0(\g)-\g\in H^{\frac{k}{\log}}$ for $k\geq 3$ and $x(\g,t)\in C([0,T];H^k)$ be a solution of (\ref{QGm},\ref{la}) with $F(x)(\g,\e,0)\in L^\infty$ and $\dg x(\g,0)\cdot \dg^2 x(\g,0)=0$. Then there exists a unique solution to (\ref{para}) with $\phi(\g,t)-\g\in C([0,T];H^{\frac{k}{\log}})$ such that $\phi(\g,0)=\phi_0(\g).$ In particular, if $\dg\phi_0(\g)>0$ then $\dg\phi(\g,t)>0$ holds for any $t\in(0,t_p]$ with $t_p\in (0,T]$.  
\end{prop}
The proof of the proposition is given in the next section. The space $H^{\frac{k}{\log}}$ is needed because we can only assume that $\lambda(x)\in H^{\frac{k}{\log}}(\T)$ for $x\in H^k$ (see the proof of Proposition \ref{proposicion}). Observe that the logarithmic modification of Sobolev norms is not a problem in the proof of the existence theorem given in \cite{FG}, because there only control of the $H^{k-1}$ norm of $\lambda(x)$ is needed, which is far away from the $H^{\frac{k}{\log}}$ norm. In the energy estimates which provide local existence, one needs to consider the integral
$$
\int \dg^k x(\g,t)\cdot \dg^k x_t(\g,t) d\g,
$$
whose most singular term coming form $\lambda(x)$ is given by
$$
I=\int \dg^k \lambda(x)(\g,t)\dg^k x(\g,t)\cdot \dg x(\g,t) d\g.
$$
Integration by parts yields
$$
I=-\int \dg^{k-1} \lambda(x)(\g,t)\dg(\dg^k x\cdot \dg x)(\g,t) d\g,
$$
and using identity \eqref{unacancelacion} one get the bound
$$
I\leq \|\dg^{k-1}\lambda(x)\|_{L^2}\|\dg(\dg^k x\cdot \dg x)\|_{L^2}\leq \|\lambda(x)\|_{H^{k-1}}\|x\|^2_{H^k}\leq C\|x\|^p_{H^k}
$$
with $p$ and $C$ constants depending on $k\geq 3$ (it is easy to observe that this extra cancellation can not be used in the $\phi$ equation). 

Next we shall show that $\xt(\xi,t)$ is a solution of \eqref{psqglz}. Here we consider $\phi$ regular enough ($\phi(\g,t)-\g\in C([0,T];H^{\frac{k}{\log}})$ with $k\geq 3$) so that it is a bona fide reparameterization satisfying \eqref{serunapara}.

The chain rule implies
\begin{equation}\label{cambio}
x_t(\g,t)=\xt_t(\phi(\g,t),t)+\phi_t(\g,t)\dxi\xt(\phi(\g,t),t).
\end{equation}
On the other hand, the equation for the evolution provides
\begin{align*}
x_t(\g,t)=&\int\frac{\dg x(\g,t)-\dg x(\e,t)}{|x(\g,t)-x(\e,t)|}d\e+\lambda(x)(\g,t)\dg x(\g,t)\\
         =&\int\frac{\dxi \xt(\phi(\g,t),t)\dg\phi(\g,t)-\dg x(\e,t)}{|x(\g,t)-x(\e,t)|}d\e+
         \lambda(x)(\g,t)\dxi \xt(\phi(\g,t),t)\dg \phi(\g,t),
 \end{align*}
and therefore       
\begin{align}
\begin{split}\label{otrac}
x_t(\g,t)=&\dxi \xt(\phi(\g,t),t)\int\frac{\dg\phi(\g,t)-\dg\phi(\e,t)}{|x(\g,t)-x(\e,t)|}d\e+\lambda(x)(\g,t)\dxi \xt(\phi(\g,t),t)\dg \phi(\g,t)\\
&+\int\frac{\dxi\xt(\phi(\g,t),t)\dg\phi(\e,t)-\dg x(\e,t)}{|x(\g,t)-x(\e,t)|}d\e.
\end{split}
\end{align}
The fact that $\phi$ is a solution of \eqref{para} together with identities (\ref{cambio},\ref{otrac}) allow us to get
$$
\xt_t(\phi(\g,t),t)=\int\frac{\dxi\xt(\phi(\g,t),t)-\dxi\xt(\phi(\e,t),t)}{|\xt(\phi(\g,t),t)-\xt(\phi(\e,t),t)|}\dg\phi(\e,t)d\e.
$$
Introducing the change of variable $\phi(\eta,t)=\zeta$ in the integral above and taking $\g=\phi^{-1}(\xi,t)$ we obtain $\xt(\xi,t)$ as a solution of \eqref{psqglz} replacing $x$ by $\xt$, $\g$ by $\xi$ and $\e$ by $\zeta$. Therefore $\tilde{x}\in C([0,T];H^{\frac{k}{\log}})$ as a consequence of Leibniz rule for derivatives of composite functions. An interesting feature in this process is the logarithm lost of derivative which affects the solutions of \eqref{psqglz}, nevertheless we will show later how to take care of that.

Once at this point one can see clearly how this reparameterization process helps to solve the following system
\begin{equation}\label{cualquiertangencial}
\xt_t(\xi,t)=\int\frac{\dxi\xt(\xi,t)-\dxi\xt(\zeta,t)}{|\xt(\xi,t)-\xt(\zeta,t)|}d\zeta+\tilde{\mu}(\xi,t)\dxi\xt(\xi,t),
\end{equation}
for any $\tilde{\mu}(\xi,t)$ having the same regularity than $\xt(\xi,t)$. We just have to repeat the argument but with the equation 
\begin{equation*}
\phi_t(\g,t)=\int_{\T}\frac{\dg\phi(\g,t)-\dg\phi(\e,t)}{|x(\g,t)-x(\e,t)|}d\e+\lambda(x)(\g,t)\dg \phi(\g,t)-\mu(\g,t),
\end{equation*}
where the function $\mu$ acts as a source term, and so long as $\phi$ and $\mu$ have the same regularity, the argument works. We then arrive to \eqref{cualquiertangencial} with $\tilde{\mu}(\xi,t)=\mu(\phi^{-1}(\xi,t),t)$. This shows that the systems \eqref{cualquiertangencial} or \eqref{psqglz} come from the system (\ref{QGm},\ref{la}) by a change of variable. 

 The main purpose of this paper is to show uniqueness for the patch problem for SQG which was until now an open problem. The following Theorem provides this result:
\begin{thm}\label{USQGp}
	Consider a solution of \eqref{weaksolution} with $\theta(x,t)$ given by a patch \eqref{patchsolution} and $D^{j}(t)$ time dependent simply connected domains whose moving boundary satisfies the arc-chord condition for any $t\in [0,T]$ and $C([0,T];C^{2,\delta}(\T))\cap C^1([0,T];C^1(\T))$ regularity. Furthermore, assume that the function $\bar{\theta}(x,t)$ given by
		\begin{equation*}
		\bar{\theta}(x,t)=\left\{\begin{array}{rl}
		\theta^1,& x\in \bar{D}^1(t),\\
		\theta^2,& x\in \bar{D}^2(t)=\R^2\setminus \bar{D}^1(t),
		\end{array}\right. 
		\end{equation*}
		satisfies \eqref{weaksolution} with $\partial\bar{D}^{j}(t)\in C([0,T];C^{2,\delta}(\T))\cap C^1([0,T];C^1(\T))$ and $\theta(x,0)=\bar{\theta}(x,0)$. Then $\theta(x,t)=\bar{\theta}(x,t)$ for any $t\in [0,T]$.
\end{thm}

This is an important part of the paper and it is proved in its section 3. In particular we show that any weak solutions of \eqref{weaksolution} given by a patch \eqref{patchsolution}, for a given parameterization \eqref{parameterization} with a certain regularity, can be reparameterized satisfying \eqref{unacancelacion}. This property is preserved in time and, together with  a new reparameterized curve, help us to  fix the tangential velocity for a contour that evolves by (\ref{QGm},\ref{la}) giving the patch solution. Then, one just needs to get uniqueness for the system (\ref{QGm},\ref{la}). Next we check the evolution of the $H^1$ Sobolev norm of the difference among two different curves evolving by (\ref{QGm},\ref{la}). We close the estimate revisiting the previous existence results and introducing new cancellation  and tools to find uniqueness by Gronwall's lemma. However in this process several different point of views respect to the previous literature are introduced.

An important linear operator in the study of patch solutions for SQG is given by
\begin{equation}\label{linealoperator}
\mathcal{L}(f)(\g)=\int_{-\pi}^\pi \frac{f(\g)- f(\g-\eta)}{|\eta|}d\eta,
\end{equation}
for $f$ $2\pi$-periodic. Since $\mathcal{L}$ is translations invariance (where we have extended $|\eta|^{-1}$ periodically), the operator is a Fourier a multiplier given by
\begin{equation}\label{lopropiedad}
\widehat{\mathcal{L}(f)}(k)=O(\log(2|k|))\widehat{f}(k),\quad\mbox{for}\quad k\in\Z\smallsetminus\{0\},\quad \widehat{\mathcal{L}(f)}(0)=0.
\end{equation}  

Uniqueness for the 2D Euler vortex patch problem was obtain in the classical Yudovich's work \cite{Yudovich}. The results presented in that paper hold in a more general setting but it is also valid for any 2D Euler weak solution with vorticity in $L^\infty(0,T;L^\infty\cap L^1)$. For the $\alpha$-system, weak solutions given by patches have been shown to be unique in \cite{KYZ}. The uniqueness result in the present paper corresponds to the more singular and physically relevant case: $\alpha=1$, but the arguments can be extended for $0<\alpha<1$. In those cases the equations for the reparameterization are more regular than \eqref{para} and there is no a logarithm derivative loss in the change of variable process. Solutions for one of the contour evolution equations were shown to be unique in \cite{FG} for $0<\al<1$. 

 Finally, in section 4 we provide uniqueness for the system \eqref{psqglz} dealing with curves without property \eqref{unacancelacion} and the evolution of the $L^2$ norm for the difference among two solutions by a different approach.

\section{Existence of an appropriate parameterization and commutator estimate}

First let us define the operators used along the proofs, namely $\dl$ and $I_{\log}$, a derivative and potential operators respectively, as the following Fourier multipliers
$$
\widehat{\dl f}(j)=\frac{j}{\log(|j|+e)}\hat{f}(j),\quad \widehat{I_{\log} f}(j)=\frac{1}{\log(|j|+e)}\hat{f}(j), 
$$ 
for $f\in L^2(\T)$. Clearly we have that $f\in L^2(\T)$ belongs to $H^{\frac{k}{\log}}$ if
$$
\dl\dg^{k-1}f\in L^2,\quad\mbox{or}\quad\dg^{k}I_{\log}f\in L^2.
$$
Next we show a commutator estimate needed in the existence and uniqueness proofs.
\begin{lemma}
Let $l^1$ be the space of absolutely convergence series. Then
\begin{equation}\label{conmutal1}
\|\dl\dg(gf)-g\dl\dg f\|_{L^2}\leq C(\|\widehat{\dg g}\|_{l^1}\|\dl f\|_{L^2}+\|\dl\dg g\|_{L^2}\|\widehat{f}\|_{l^1}),
\end{equation} where $C$ is a universal constant. In particular Sobolev's embedding implies that for any $\epsilon>0$ there is a constant $C_\epsilon>0$ such that
\begin{equation}\label{conmutador}
\|\dl\dg(gf)-g\dl\dg f\|_{L^2}\leq C_\ep(\|g\|_{H^{3/2+\epsilon}}\|\dl f\|_{L^2}+\|\dl\dg g\|_{L^2}\|f\|_{H^{1/2+\epsilon}}).
\end{equation} 
\end{lemma}
\emph{Proof:}
We have that 
$$
|(\dl\dg(gf)-g\dl\dg f)\hat{\,}(j)|\leq \sum_{l}\Big|\frac{j^2}{\log(|j|+e)}-\frac{(j-l)^2}{\log(|j-l|+e)}\Big||\hat{f}(j-l)||\hat{g}(l)|,
$$
and the function $h(j)=j^2/\log(|j|+e)$ satisfies
$$
h(j)-h(j-l)=\int_0^1 \frac{d}{dr}h((j-l)+rl)dr=l\int_0^1h'(rl+(j-l))dr, 
$$  
and therefore
$$
|h(j)-h(j-l)|\leq |l|\int_0^1|h'(rl+(j-l))|dr\leq \frac{3(|l|+|j-l|)}{\log(|l|+|j-l|+e)}|l|.
$$
It yields
$$
|h(j)-h(j-l)|\leq\frac{3|l|^2}{\log(|l|+e)}+\frac{3|l||j-l|}{\log(|j-l|+e)},
$$
and finally
$$
|(\dl\dg(gf)-g\dl\dg f)\hat{\,}(j)|\leq \sum_{l}\frac{3|l|^2}{\log(|l|+e)}|\hat{f}(j-l)||\hat{g}(l)|+
\sum_{l}\frac{3|l||j-l|}{\log(|j-l|+e)}|\hat{f}(j-l)||\hat{g}(l)|.
$$
Then Parseval's Theorem gives
\begin{align*}
\|\dl\dg(gf)-g\dl\dg f\|_{L^2}\leq & \Big(\sum_{j}\Big(\sum_{l}\frac{3|l|^2}{\log(|l|+e)}|\hat{f}(j-l)||\hat{g}(l)|\Big)^2\Big)^{1/2}\\
&+\Big(\sum_j\Big(\sum_{l}\frac{3|l||j-l|}{\log(|j-l|+e)}|\hat{f}(j-l)||\hat{g}(l)|\Big)^2\Big)^{1/2}.
\end{align*}
Minkowski inequality provides \eqref{conmutal1}. The proof ends by Sobolev's embedding in dimension one.\\
\\
\emph{Proof of Proposition \ref{proposicion}:}
Without lost of generality we may consider the case $k=3$, because the extension to $k>3$ is just a straightforward exercise once we know how to handle $k=3$. Also, in order to be concise we will show only the main part of the proof. That is, we will deal with the more dangerous terms in the needed estimates, leaving as an exercise to the reader the treatment to all the others more benevolent characters. 
In the main core of the proof are energy estimates, from them and with nowadays well-known  mollifying arguments one can apply the classical Picard to conclude existence. The whole strategy can be found in \cite{bertozzi-Majda}, chapter 3. 

Often, in the following we will have to write double integrals in variables, say $\g$ and $\eta$, and differences $f(\g)-f(\g-\eta)$. To simplify notation we shall write $f=f(\g,t)$, $f'=f(\gamma-\eta,t)$ and $f-f'=f_-$ when there is no danger of confusion. Furthermore, we shall write $\int=\int_{\T}$ and $id$ is the identity, $C(t)$ will be a polynomial function in $\|F(x)\|_{L^\infty}$ and $\|x\|_{H^3}$ so that $C(t)\in C([0,T])$. As was mentioned before, most of the time we will show how to estimate the most singular terms: those in which the derivative of higher order is involved by the use of the Leibnitz's derivative rule. The rest of the terms are denoted by $l.o.t.$ standing for lower order terms. Writing $l.o.t.\in X$ means that the lower order terms belong to the space $X$.    

First we consider the evolution of the $L^2$ norm:
\begin{align*}
\frac12\frac{d}{dt}\|\phi-id\|_{L^2}^2=\int (\phi-id)\phi_t d\g=I_1+I_2,
\end{align*} 
where
$$
I_1=\int (\phi-id)\int\frac{\dg\phi_-}{|x_-|}d\e d\g,\quad I_2=\int \lambda(x)(\phi-id) \dg \phi d\g.
$$
For $I_1$ we find
\begin{align*}
I_1=&\int\!\!\!\int (\phi(\g,t)\!-\!\g)\frac{\dg\phi(\g,t)\!-\!\de\phi(\e,t)}{|x(\g,t)\!-\!x(\e,t)|}d\e d\g=-\!\int\!\!\!\int (\phi(\e,t)\!-\!\e)\frac{\dg\phi(\g,t)\!-\!\de\phi(\e,t)}{|x(\g,t)\!-\!x(\e,t)|}d\e d\g\\
=&\frac12\int\!\!\!\int (\phi(\g,t)-\g-(\phi(\e,t)-\e))\frac{\dg\phi(\g,t)-\de\phi(\e,t)}{|x(\g,t)-x(\e,t)|}d\e d\g,
\end{align*}
hence
\begin{align*}
I_1=&\frac12\int\!\!\!\int (\phi-id)_-\frac{\dg((\phi-id)_-)}{|x_-|}d\e d\g.
\end{align*}
Integration by parts yields
\begin{align*}
I_1=&\frac14\int\!\!\!\int |(\phi-id)_-|^2 \frac{x_-\cdot\dg x_-}{|x_-|^3}d\e d\g.
\end{align*}
Now we use \eqref{unacancelacion} to rewrite
$$
\frac{x_-\cdot\dg x_-}{|x_-|^3}=\frac{x_-\cdot\dg x_--\dg x\cdot\dg^2 x\e^2}{|x_-|^3},
$$
and obtain
\begin{equation}\label{cancelsuper}
\Big|\frac{x_-\cdot\dg x_-}{|x_-|^3}\Big|\leq \frac{2\|x\|^2_{C^{2,\frac12}}|\e|^{2+\frac12}}{|
x_-|^3}\leq2\|F(x)\|^3_{L^\infty}\|x\|^2_{C^{2,\frac12}}|\e|^{-\frac12}.
\end{equation}
This yields
\begin{align*}
I_1\leq&\frac12\|F(x)\|^3_{L^\infty}\|x\|^2_{C^{2,\frac12}}\int|\e|^{-\frac12}\int |(\phi-id)_-|^2  d\g d\e\\
\leq&\|F(x)\|^3_{L^\infty}\|x\|^2_{C^{2,\frac12}}\int|\e|^{-\frac12}\int (|\phi-id|^2+|(\phi-id)'|^2)d\g d\e\\
\leq&2\|F(x)\|^3_{L^\infty}\|x\|^2_{C^{2,\frac12}}\|\phi-id\|^2_{L^2}(t)\leq C(t)\|\phi-id\|^2_{L^2}(t).
\end{align*}
 The term $I_2$ can be rewritten as follows
$$
I_2=\int \lambda(x)(\phi-id) (\dg \phi-1)d\g+\int \lambda(x)(\phi-id)d\g.
$$
The first term above can be handled by integration by parts. In the second Cauchy-Schwarz inequality yields
$$
I_2\leq C(t)\|\phi-id\|^2_{L^2}(t)+\frac12\|\lambda(x)\|^2_{L^2}.
$$
The bounds for $\lambda(x)$ (below we show that $\lambda\in H^{\frac{3}{\log}}$) finally provide
\begin{align}\label{phiL2}
\frac{d}{dt}\|\phi-id\|^2_{L^2}(t)\leq C(t)\|\phi-id\|^2_{L^2}(t)+C(t).
\end{align}
Next, we consider the evolution of the higher order norm 
\begin{align}
\begin{split}\label{JK}
\frac12\frac{d}{dt}\|\dl\dg^2\phi\|_{L^2}^2&=  \int\dl\dg^2\phi\dl\dg^2\phi_t d\g\\
=\int\dl\dg^2\phi&\dl\Big(\dg^2\Big(\int\frac{\dg \phi_-}{|x_-|}d\eta\Big)\Big)d\g+\int\dl\dg^2\phi\dl\Big(\dg^2\Big(\lambda(x)\dg\phi\Big)\Big)d\g\\
=J+K&,
\end{split}
\end{align}to bound the $J$ and $K$ terms. 

With $J$ we split further $J=J_1+J_2+J_3$ where 
$$
J_1=\int\dl\dg^2\phi\dl\Big(\int\frac{\dg^3 \phi_-}{|x_-|}d\eta\Big)d\g,\quad J_2=-2\int\dl\dg^2\phi\dl\Big(\int\frac{\dg^2 \phi_-x_-\cdot\dg x_-}{|x_-|^3}d\eta\Big)d\g,
$$
and
$$J_3=-\int\dl\dg^2\phi\dl\Big(\int\dg \phi_-\dg\Big(\frac{x_-\cdot\dg x_-}{|x_-|^3}\Big)d\eta\Big)d\g.$$
The fact that $|\dg x| $ does not depend on $\gamma$ gives 
$$
\int\dl\dg^2\phi\dl\Big(\int  \frac{\dg^3 \phi_-}{|\dg x||\eta|}\Big)d\eta\Big)d\g=
\frac1{|\dg x|}\int\dl\dg^2\phi\dl \mathcal{L}(\dg^3\phi)d\g=0
$$
where $\mathcal{L}$ was defined in \eqref{linealoperator} and has properties \eqref{lopropiedad}. Therefore one obtains 
$$
J_1=\int\dl\dg^2\phi\dl\Big(\int \dg^3 \phi_- \Big(\frac{1}{|x_-|}-\frac{1}{|\dg x||\eta|}\Big)d\eta\Big)d\g.
$$
This extra cancellation suggest the further splitting  $J_1=J_{1,1}+J_{1,2}$ where 
$$
J_{1,1}=\int\dl\dg^2\phi\dl\Big(\dg^3\phi\int\Big(\frac{1}{|x_-|}-\frac{1}{|\dg x||\eta|}\Big)d\eta\Big)d\g,
$$
$$
J_{1,2}=-\int\dl\dg^2\phi\dl\Big(\int\dg^3\phi'\Big(\frac{1}{|x_-|}-\frac{1}{|\dg x||\eta|}\Big)d\eta\Big)d\g,
$$
and $J_{1,1}=J_{1,1}^1+J_{1,1}^2+J_{1,1}^3$ where
$$
J_{1,1}^1=\int\dl\dg^2\phi[\dl\dg(A\dg^2\phi)-A\dl\dg(\dg^2\phi)]d\g,
$$
$$
J_{1,1}^2=\int\dl\dg^2\phi A\dl\dg(\dg^2\phi)d\g,\quad J_{1,1}^3=-\int\dl\dg^2\phi \dl(\dg A\dg^2\phi)d\g,
$$
with 
\begin{equation}\label{Adef}
A=\int\Big(\frac{1}{|x_-|}-\frac{1}{|\dg x||\eta|}\Big)d\eta.
\end{equation}

In $J_{1,1}^1$ we use the commutator estimate \eqref{conmutador} to find
$$
J_{1,1}^1\leq  C\|\dl\dg^2\phi \|_{L^2}(\|A\|_{H^{2}}\|\dl\dg^2\phi \|_{L^2}+\|\dl\dg A\|_{L^2}\|\dg^2\phi\|_{H^{\frac{1}{\log}}})\leq 
C\|A\|_{H^{2}}\|\dl\dg^2\phi \|^2_{L^2}.
$$ 
Furthermore we have 
\begin{equation}\label{Ad}
\dg A=- \int\frac{x_-\cdot \dg x_--\dg x\cdot\dg^2 x\eta^2}{|x_-|^3}d\eta,
\end{equation}
and therefore
$$
\dg^2 A=-\int\frac{x_-\cdot \dg^2 x_--\dg x\cdot\dg^3 x\eta^2}{|x_-|^3}d\eta+l.o.t.
$$
where $\|l.o.t.\|_{L^2}\leq C(t)$. Identity \eqref{unacancelacion} yields
\begin{align*}
x_-\cdot \dg^2 x_--\dg x\cdot\dg^3 x\eta^2&=(x_--\dg x\eta)\cdot \dg^2 x_-+\eta\dg x\cdot \dg^2 x_--\dg x\cdot\dg^3 x\eta^2\\
&=(x_--\dg x\eta)\cdot \dg^2 x_--\eta\dg x_-\cdot \dg^2 x'+
|\dg^2 x|^2\eta^2,
\end{align*}
implying
$$
x_-\cdot \dg^2 x_--\dg x\cdot\dg^3 x\eta^2=(x_--\dg x\eta)\cdot \dg^2 x_--\eta(\dg x_--\dg^2x\eta)\cdot \dg^2 x'+\eta^2\dg^2x\cdot\dg^2 x_-.
$$
The above configuration provides
$$
|x_-\cdot \dg^2 x_--\dg x\cdot\dg^3 x\eta^2|\leq 3\|x\|_{C^2}\|x\|_{C^{2,\frac12}}|\eta|^{2+\frac12},
$$
and therefore
$$
\Big|-\int\frac{(x_--\dg x\eta)\cdot \dg^2 x_--\dg x\cdot\dg^3 x\eta^2}{|x_-|^3}d\eta\Big|\leq 3\|x\|_{C^2}\|x\|_{C^{2,\frac12}}\|F(x)\|^3_{L^\infty},
$$
implying that $\dg^2 A\in C([0,T],L^2)$ and the estimate
$$
J_{1,1}^1\leq C(t)\|\dl\dg^2\phi \|^2_{L^2}.
$$
Then, integration by parts yields
$$
J_{1,1}^2\leq \|\dg A\|_{L^\infty}\|\dl\dg^2\phi \|^2_{L^2}\leq \|A\|_{H^2}\|\dl\dg^2\phi \|^2_{L^2}\leq  C(t)\|\dl\dg^2\phi \|^2_{L^2}.
$$
In order to estimate $J_{1,1}^3$ we use the following inequalities
\begin{align*}
J_{1,1}^3&\leq \|\dl\dg^2\phi \|_{L^2}\|\dl(\dg A\dg^2\phi)\|_{L^2}\leq C\|\dl(\dg A)\|_{L^2}\|\dl\dg^2\phi \|^2_{L^2}\leq C\|A\|_{H^2}\|\dl\dg^2\phi \|^2_{L^2}\\
&\leq C(t)\|\dl\dg^2\phi \|^2_{L^2}.
\end{align*}
Hence
$$
J_{1,1}\leq C(t)\|\dl\dg^2\phi \|^2_{L^2}.
$$
It remains to control $J_{1,2}$. We rewrite $J_{1,2}=J_{1,2}^1+J_{1,2}^2$ given by
$$
J_{1,2}^1=\int\dl\dg^2\phi\dl\Big(\int\dg^3\phi'\Omega_1d\eta\Big)d\g,\quad 
J_{1,2}^2=\int\dl\dg^2\phi\dl\Big(\int\dg^3\phi'\Omega_2d\eta\Big)d\g
$$
where $|\dg x|^{-1}|\eta|^{-1}-|x_-|^{-1}=\Omega_1+\Omega_2$ with
\begin{equation}\label{Omega1y2}
\Omega_1=\frac{\big|\frac{x_-}{\eta}-\dg x\big|^2}{|\dg x|\big(|\dg x|+|\frac{x_-}{\eta}|\big)|x_-|},\quad \Omega_2=2\frac{\big(\frac{x_-}{\eta}-\dg x\big)\cdot\dg x}{|\dg x|\big(|\dg x|+|\frac{x_-}{\eta}|\big)|x_-|}.
\end{equation}
Next we will show how to deal with $J_{1,2}^2$ and since the kernel $\Omega_2$ is more singular than $\Omega_1$, we leave to the reader the analogous details for $J_{1,2}^1$. 

Identity \eqref{unacancelacion} allows us to rewrite
\begin{equation}\label{omega2}
\Omega_2=2\frac{\big(\frac{x_-}{\eta}-\dg x+\frac12\dg^2 x\eta\big)\cdot\dg x}{|\dg x|\big(|\dg x|+|\frac{x_-}{\eta}|\big)|x_-|}.
\end{equation}
and the splitting $J_{1,2}^2=J_{1,2}^{2,1}+J_{1,2}^{2,2}$ where
\begin{equation}\label{J12s2d}
J_{1,2}^{2,1}=\int\dl\dg^2\phi\dl\Big(\int\dg^2\phi'\de\Omega_2d\eta\Big)d\g\quad 
J_{1,2}^{2,2}=-\int\dl\dg^2\phi\dl\Big(\dg^2\phi'\Omega_2\Big|_{\eta=-\pi}^{\eta=\pi}\Big)d\g.
\end{equation}
In the case of $J_{1,2}^{2,2}$ let us observe that the functions $\Omega_2(\g,\pm\pi)$ are regular enough to obtain 
\begin{align*}
J_{1,2}^{2,2}&\leq \|\dl\dg^2\phi \|_{L^2}\|\dl(\dg^2\phi'\Omega_2\Big|_{\eta=-\pi}^{\eta=\pi})\|_{L^2}\leq C\|\dl\dg^2\phi\|^2_{L^2}
\sum_{\pm}\|\dl\Omega_2\,_{|_{\eta=\pm\pi}}\|_{L^2}\\
&\leq C\|\dl\dg^2\phi\|^2_{L^2}\sum_{\pm}
\|\dg\Omega_2\,_{|_{\eta=\pm\pi}}\|_{L^2}\leq C(t)\|\dl\dg^2\phi\|^2_{L^2}.
\end{align*}
Regarding to $J_{1,2}^{2,1}$, we proceed as follows
\begin{align}
\begin{split}\label{J12s21cota}
J_{1,2}^{2,1}&\leq \|\dl\dg^2\phi \|_{L^2}\|\dl\Big(\int\dg^2\phi'\de\Omega_2d\eta\Big)\|_{L^2}\leq \|\dl\dg^2\phi\|_{L^2}\|\dg\Big(\int\dg^2\phi'\de\Omega_2d\eta\Big)\|_{L^2}\\
&\leq \|\dl\dg^2\phi\|_{L^2}\Big(\|\int\dg^2\phi'\dg\de\Omega_2d\eta\|_{L^2}+\|\int\dg^3\phi'\de\Omega_2d\eta\|_{L^2}\Big).
\end{split}
\end{align}
And two new terms appear that have to be controlled in $L^2$:
\begin{equation}\label{BD}
B=\int\dg^2\phi'\dg\de\Omega_2d\eta,\quad D=\int\dg^3\phi'\de\Omega_2d\eta.
\end{equation}
In order to do that first we will prove the bound $\|\dg\de\Omega_2\|_{L^2}\leq C(t)$ to obtain
\begin{equation}\label{Bcota}
\|B\|_{L^2}\leq \|\dg^2\phi\|_{L^1}\|\dg\de\Omega_2\|_{L^2}\leq C(t)\|\dl\dg^2\phi\|_{L^2}.
\end{equation}
With the help of formula \eqref{omega2} we split $\de\Omega_2=\p\Omega_{2,1}+\p\Omega_{2,2}+\p\Omega_{2,3}+\p\Omega_{2,4}$
where
$$
\p\Omega_{2,1}=\frac{-2}{\eta^2}\frac{(x_--\dg x\eta+\frac12\dg^2 x\eta^2)\cdot\dg x}{|\dg x|\big(|\dg x|+|\frac{x_-}{\eta}|\big)|x_-|},\quad 
\p\Omega_{2,2}=\frac{2}{\eta}\frac{(\dg x'-\dg x+\dg^2 x\eta)\cdot\dg x}{|\dg x|\big(|\dg x|+|\frac{x_-}{\eta}|\big)|x_-|},
$$
\begin{equation}\label{Omega21234}
\p\Omega_{2,3}=\frac{-2}{\eta}\frac{(x_--\dg x\eta+\frac12\dg^2 x\eta^2)\cdot\dg x}{|\dg x|\big(|\dg x|+|\frac{x_-}{\eta}|\big)|x_-|^2}\frac{x_-\cdot\dg x'}{|x_-|},
\end{equation}
and
$$ 
\p\Omega_{2,4}=\frac{-2}{\eta}\frac{(x_--\dg x\eta+\frac12\dg^2 x\eta^2)\cdot\dg x}{|\dg x|\big(|\dg x|+|\frac{x_-}{\eta}|\big)^2|x_-|}\D\frac{\D\frac{x_-}{\eta}\cdot\frac{\dg x'\eta-x_-}{\eta^2}}{\D\big|\frac{x_-}{\eta}\big|}.
$$
Next we will show how to deal with $\dg\p\Omega_{2,1}$ and since the other kernels are similar or even easier to handle we will skip the details. 

We have
$$
\dg\p\Omega_{2,1}=\frac{-2}{\eta^2}\frac{(\dg x_--\dg^2 x\eta+\frac12\dg^3 x\eta^2)\cdot\dg x}{|\dg x|\big(|\dg x|+|\frac{x_-}{\eta}|\big)|x_-|}+l.o.t.,
$$
where $\|l.o.t.\|_{L^2}\leq C(t)$. The identity
$$
\dg x_--\dg^2 x\eta+\frac12\dg^3 x\eta^2=\eta^2\int_0^1r(\dg^3x(\g)-\dg^3x(\g+(r-1)\eta))dr
$$
allows us to write
\begin{align*}
(\dg x_--\dg^2 x\eta+\frac12\dg^3 x\eta^2)\cdot\dg x=&\eta^2\int_0^1r \dg^3x(\g\!+\!(r\!-\!1)\eta))\cdot(\dg x(\g\!+\!(r\!-\!1)\eta)-\dg x(\g))dr\\
+
\eta^2\int_0^1r&(\dg^3x(\g)\cdot\dg x(\g)-\dg^3x(\g\!+\!(r\!-\!1)\eta)\cdot\dg x(\g\!+\!(r\!-\!1)\eta))dr.
\end{align*}
The use of equality \eqref{unacancelacion} and integration by parts in $r$ yield
\begin{align*}
(\dg x_--\dg^2 x\eta+\frac12\dg^3 x\eta^2)\cdot\dg x=&\eta^2\int_0^1r \dg^3x(\g\!+\!(r\!-\!1)\eta))\cdot(\dg x(\g\!+\!(r\!-\!1)\eta)-\dg x(\g))dr\\
&-
\eta^3\int_0^1r^2\dg^2x(\g\!+\!(r\!-\!1)\eta)\cdot\dg^3 x(\g\!+\!(r\!-\!1)\eta)dr,
\end{align*}
and therefore
\begin{align*}
|(\dg x_--\dg^2 x\eta+\frac12\dg^3 x\eta^2)\cdot\dg x|\leq &2|\eta|^3\|x\|_{C^2}\int_0^1 |\dg^3x(\g\!+\!(r\!-\!1)\eta)|dr.
\end{align*}
Hence
$$
|\dg\p\Omega_{2,1}|\leq 4\|F(x)\|_{L^\infty}^3\|x\|_{C^2}\int_0^1 |\dg^3x(\g\!+\!(r\!-\!1)\eta)|dr+|l.o.t.|.
$$
Finally an integration in $\gamma$ gives the desired property:  $\|\dg\p\Omega_{2,1}\|_{L^2}\leq C(t).$ Analogously we have $\|\dg\p\Omega_{2,j}\|_{L^2}\leq C(t)$ for $j=2,3,4$ and therefore the same bound holds for $\dg\de\Omega_2$:
$$
\|\dg\de\Omega_{2}\|_{L^2}\leq C(t).
$$
We achieve the desired estimate \eqref{Bcota}. 

 Regarding $D$, we first integrate by parts and then split
$$
D=\int\de(\dg^2\phi_-)\de\Omega_2d\eta=-\int\dg^2\phi_-\de^2\Omega_2d\eta+\dg^2\phi_-\de\Omega_2\Big|_{\eta=-\pi}^{\eta=\pi}=D_1+D_2.
$$
Then formulas \eqref{Omega21234} show that the functions $\de\Omega_2(\g,\pm\pi)$ are regular enough to get an appropriate bound for $D_2$: 
$$
\|D_2\|_{L^2}\leq C(t)\|\dg^2\phi\|_{L^2}\leq  C(t)\|\dl\dg^2\phi\|_{L^2}.
$$ 
Following the decomposition for $\de \Omega_2$ in \eqref{Omega21234}, let us introduce $\de\p\Omega_{2,1}=\p^2\Omega_{2,1}^1+\p^2\Omega_{2,1}^2+\p^2\Omega_{2,1}^3+\p^2\Omega_{2,1}^4$ 
where
$$
\p^2\Omega_{2,1}^1=\frac{4}{\eta^3}\frac{(x_--\dg x\eta+\frac12\dg^2 x\eta^2)\cdot\dg x}{|\dg x|\big(|\dg x|+|\frac{x_-}{\eta}|\big)|x_-|},\quad 
\p^2\Omega_{2,1}^2=\frac{-2}{\eta^2}\frac{(-\dg x_-+\dg^2 x\eta)\cdot\dg x}{|\dg x|\big(|\dg x|+|\frac{x_-}{\eta}|\big)|x_-|},
$$
$$
\p^2\Omega_{2,1}^3=\frac{2}{\eta^2}\frac{(x_--\dg x\eta+\frac12\dg^2 x\eta^2)\cdot\dg x}{|\dg x|\big(|\dg x|+|\frac{x_-}{\eta}|\big)|x_-|^2}\frac{x_-\cdot\dg x'}{|x_-|},
$$
and
$$
\p^2\Omega_{2,1}^4=\frac{2}{\eta^2}\frac{(x_--\dg x\eta+\frac12\dg^2 x\eta^2)\cdot\dg x}{|\dg x|\big(|\dg x|+|\frac{x_-}{\eta}|\big)^2|x_-|}\D\frac{\D\frac{x_-}{\eta}\cdot\frac{\dg x'\eta-x_-}{\eta^2}}{\D\big|\frac{x_-}{\eta}\big|}.
$$
As was shown before, we have
$$
x_--\dg x\eta+\frac12\dg^2 x\eta^2=\frac12\eta^3\int_0^1r^2\dg^3x(\g+(r-1)\eta)dr
$$
and therefore
$$
|\p^2\Omega_{2,1}^1|\leq \frac2{|\eta|}\|F(x)\|_{L^\infty}^2\int_0^1|\dg^3x(\g+(r-1)\eta)|dr.
$$
Analogously, we obtain
$$
\sum_{j=2}^4|\p^2\Omega_{2,1}^j|\leq \frac2{|\eta|}(\|F(x)\|_{L^\infty}^2\!+\!\|F(x)\|_{L^\infty}^3\|x\|_{C^1}\!+\!\|F(x)\|_{L^\infty}^2\|x\|_{C^1})\!\int_0^1\!|\dg^3x(\g\!+\!(r\!-\!1)\eta)|dr,
$$
implying
$$
|\de\p\Omega_{2,1}|\leq \frac1{|\eta|}C(t)\int_0^1|\dg^3x(\g+(r-1)\eta)|dr.
$$
The same approach for $\de\p\Omega_{2,j}$ with $j=2,3,4$ yields
$$
|\de^2\Omega_{2}|\leq \frac1{|\eta|}C(t)\int_0^1|\dg^3x(\g+(r-1)\eta)|dr.
$$
Therefore we get the estimate
$$
|D_1|\leq C(t)\|\dg^2\phi\|_{C^{\frac13}}\int\frac1{|\eta|^{2/3}} \int_0^1|\dg^3x(\g+(r-1)\eta)|drd\eta,
$$ 
and consequently  
$$
\|D_1\|_{L^2}\leq C(t)\|\dg^2\phi\|_{C^{\frac13}}\|\dg^3x\|_{L^2}\leq C(t)\|\dg^2\phi\|_{H^{\frac{11}{12}}}\leq 
C(t)\|\dl\dg^2\phi\|_{L^2}
$$
by Sobolev embedding. Putting all those estimates together we obtain
$$
\|D\|_{L^2}\leq C(t)\|\dl\dg^2\phi\|_{L^2},
$$
which together with \eqref{Bcota} allows us to get finally the needed estimate for $J_{1,2}^{2,1}$ in \eqref{J12s2d} using \eqref{J12s21cota}. We are then done with $J_{1,2}^2$. 

For the less singular kernel $\Omega_1$ in \eqref{Omega1y2} a similar analysis yields
$$
J_{1,2}^1\leq C(t)\|\dl\dg^2\phi\|^2_{L^2}.
$$
Hence the same estimate is achieved for $J_{1,2}$ and accordingly for $J_1$:
$$
J_1\leq C(t)\|\dl\dg^2\phi\|^2_{L^2}.
$$

Next we estimate $J_2=J_{2,1}+J_{2,2}$ given by
$$
J_{2,1}=-2\int\dl\dg^2\phi\dl\Big(\dg^2 \phi\int\frac{x_-\cdot\dg x_-}{|x_-|^3}d\eta\Big)d\g=2\int\dl\dg^2\phi\dl\Big(\dg^2\phi\dg A\Big)d\g,
$$
and
$$
J_{2,2}=2\int\dl\dg^2\phi\dl\Big(\int\dg^2\phi'\frac{x_-\cdot\dg x_-}{|x_-|^3}d\eta\Big)d\g=
2\int\dl\dg^2\phi\dl\Big(\int\dg^2\phi'\Omega_3d\eta\Big)d\g,
$$
where $\dg A$ was introduced in \eqref{Ad} and the kernel $\Omega_3$ can be rewritten as
$$
\Omega_3=\frac{x_-\cdot\dg x_--\dg x\cdot\dg^2 x\eta^2}{|x_-|^3}.
$$ Observe that $J_{2,1}=2J_{1,1}^3$ and therefore we already know the estimate of that term. The other $J_{2,2}$ is similar to $J_{1,2}^{2,1}$ because the kernel $\Omega_3$ is of degree $0$ as $\partial_\eta\Omega_2$, and has the same lost of regularity in the tangential direction. Then, as before we obtain
$$
|\dg \Omega_3|\leq C(t)\int_0^1|\dg^3 x(\g+(r-1)\eta)|dr,\quad |\de \Omega_3|\leq \frac1{|\eta|}C(t)\int_0^1|\dg^3 x(\g+(r-1)\eta)|dr.
$$
helping to estimate $J_{2,2}$, and
$$
J_2\leq C(t)\|\dl\dg^2\phi\|^2_{L^2}.
$$
Finally, to deal with $J_3$, we proceed as follows
\begin{align*}
J_3&
\leq \|\dl\dg^2\phi\|_{L^2}\|\dl\Big(\int\dg \phi_-\dg\Big(\frac{x_-\cdot\dg x_-}{|x_-|^3}\Big)d\eta\Big)\|_{L^2}\\
&\leq \|\dl\dg^2\phi\|_{L^2}\|\dg\Big(\int\dg \phi_-\dg\Big(\frac{x_-\cdot\dg x_-}{|x_-|^3}\Big)d\eta\Big)\|_{L^2},
\end{align*}
that is
$$
J_3\leq \|\dl\dg^2\phi\|_{L^2}\Big(\|\int\dg^2 \phi_-\dg\Big(\frac{x_-\cdot\dg x_-}{|x_-|^3}\Big)d\eta \|_{L^2}
+\|\int\dg \phi_-\dg^2\Big(\frac{x_-\cdot\dg x_-}{|x_-|^3}\Big)d\eta \|_{L^2}\Big).
$$
Next let us observe that the two inequalities
$$
\Big|\dg^2 \phi_-\dg\Big(\frac{x_-\!\cdot\!\dg x_-}{|x_-|^3}\Big)\Big|\leq \frac{\|\dg^2\phi\|_{C^{\frac13}}}{|\eta|^{\frac23}}\Big(\|F(x)\|_{L^\infty}^2\!\int_0^1\!\!|\dg^3 x(\g\!+\!(r\!-\!1)\eta)|dr\!+\!4\|F(x)\|_{L^\infty}^3\|x\|_{C^2}^2\Big),
$$
$$
\Big|\dg \phi_-\dg^2\Big(\frac{x_-\cdot\dg x_-}{|x_-|^3}\Big)\Big|\leq C(t)\|\dg^2\phi\|_{L^\infty}\Big(|\dg^3 x|+|\dg^3 x'|+\int_0^1\!\!|\dg^3 x(\g\!+\!(r\!-\!1)\eta)|dr+1\Big),
$$
together with Sobolev embedding yield
$$
J_3\leq C(t)\|\dl\dg^2\phi\|_{L^2}(\|\dg^2\phi\|_{C^{\frac13}}+\|\dg^2\phi\|_{L^\infty})\leq C(t)\|\dl\dg^2\phi\|_{L^2}^2,
$$
giving us the control:
$$
J\leq C(t)\|\dl\dg^2\phi\|_{L^2}^2.
$$

To finish, it remains to deal with $K$. First we will show the regularity of $\lambda(x)\in C([0,T];H^{\frac{3}{\log}})$. To do that we begin observing that $\lambda(x)\in C([0,T];L^2)$, and continue showing that $I_{\log}(\dg^3\lambda(x))\in C([0,T];L^2)$ with the following decomposition $\dg^3\lambda(x)=E_1+E_2+E_3$ where
$$
E_1=-\frac{\dg^3x}{|\dg x|^2}\cdot\dg\Big(\int\frac{\dg x_-}{|x_-|}d\eta\Big),\quad
E_2=-2\frac{\dg^2x}{|\dg x|^2}\cdot\dg^2\Big(\int\frac{\dg x_-}{|x_-|}d\eta\Big),
$$
and
$$
E_3=-\frac{\dg x}{|\dg x|^2}\cdot\dg^3\Big(\int\frac{\dg x_-}{|x_-|}d\eta\Big).
$$
The inequality
$$
E_1\leq (\|F(x)\|_{L^\infty}^3\|\dg^2x\|_{C^{\frac12}}+\|F(x)\|_{L^\infty}^4\|\dg^2x\|^2_{L^\infty})|\dg^3x|
$$
gives $E_1\in C([0,T];L^2)$. For $E_2$ we consider
$E_2=E_{2,1}+E_{2,2}+E_{2,3}$ where
$$
E_{2,1}=-2\frac{\dg^2x}{|\dg x|^2}\cdot\Big(\int\frac{\dg^3 x_-}{|x_-|}d\eta\Big),\quad 
E_{2,2}=-2\frac{\dg^2x}{|\dg x|^2}\cdot\Big(\int\frac{\dg^2 x_-}{|x_-|^3}x_-\cdot\dg x_-d\eta\Big)
$$
and
$$
E_{2,3}=-2\frac{\dg^2x}{|\dg x|^2}\cdot\Big(\int \dg^2 x_-\dg\Big(\frac{x_-\cdot\dg x_-}{|x_-|^3}\Big)d\eta\Big).
$$
A similar approach provides $E_{2,2}$ and $E_{2,3}$ in $C([0,T];L^2)$.As usual we will focus our attention in the most singular 
term $E_{2,1}$, which can be decomposed as $E_{2,1}=E_{2,1}^1+E_{2,1}^2+E_{2,1}^3$ where
$$
E_{2,1}^1=-2\frac{\dg^2x}{|\dg x|^2}\cdot\Big(\int\dg^3 x_-\Big(\frac{1}{|x_-|}-\frac{1}{|\dg x||\eta|}\Big)d\eta\Big),\quad 
E_{2,1}^2=\frac{2}{|\dg x|^3}\Big(\int\frac{\dg^2x_-\cdot\dg^3 x'}{|\eta|}d\eta\Big),
$$
and
\begin{equation}\label{e213}
E_{2,1}^3=\frac{-2}{|\dg x|^3}\int \frac{(\dg^2x\cdot\dg x^3)_-}{|\eta|}d\eta=\frac{-2}{|\dg x|^3}\mathcal{L}(\dg^2x\cdot\dg^3 x).
\end{equation}
As before one finds
$$
|E_{2,1}^1|+|E_{2,1}^2|\leq C(\|F(x)\|_{L^\infty}^4\|x\|_{C^{2,\frac12}}^2+1)(|\dg^3 x|+\int |\eta|^{-\frac12}|\dg^3 x'|d\eta),
$$
and consequently $E_{2,1}^1, E_{2,1}^2\in C([0,T];L^2)$. It remains then to deal with $E_{2,1}^3$, which is the most singular term not belonging to $C([0,T];L^2)$. Nevertheless one has
$$
\|I_{\log}(E_{2,1}^3)\|_{L^2}\leq 2\|F(x)\|_{L^\infty}^3\|I_{\log}(\mathcal{L}(\dg^2x\cdot\dg^3 x))\|_{L^2}\leq 
C\|F(x)\|_{L^\infty}^3\|\dg^2x\cdot\dg^3 x\|_{L^2}
$$
as a consequence of properties \eqref{lopropiedad}, from where we reach the desired estimate
$$
\|I_{\log}(E_{2,1}^3)\|_{L^2}\leq C\|F(x)\|_{L^\infty}^3\|\dg^2x\|_{L^\infty}\|\dg^3 x\|_{L^2}\leq C(t).
$$

In the following, we show that all the remaining terms (except one) are integrable in $C([0,T];L^2)$. This singular term is a constant times $E_{2,1}^3$. We are done with $E_{2,1}$ and consequently with $E_2$. 

Regarding $E_{3}$, we introduce the splitting $E_3=E_{3,1}+E_{3,2}+E_{3,3}+E_{3,4}$ where
$$
E_{3,1}=-\frac{\dg x}{|\dg x|^2}\cdot\int\frac{\dg^4 x_-}{|x_-|}d\eta,\quad 
E_{3,2}=3\frac{\dg x}{|\dg x|^2}\cdot\int\frac{\dg^3 x_-}{|x_-|^3}x_-\cdot\dg x_-d\eta,
$$
and
$$
E_{3,3}=3\frac{\dg x}{|\dg x|^2}\cdot\int\dg^2 x_-\dg\Big(\frac{x_-\cdot\dg x_-}{|x_-|^3}\Big)d\eta,\quad 
E_{3,4}=\frac{\dg x}{|\dg x|^2}\cdot\int\dg x_-\dg^2\Big(\frac{x_-\cdot\dg x_-}{|x_-|^3}\Big)d\eta.
$$
Using \eqref{cancelsuper}, $E_{3,2}$ has the following estimate
$$
|E_{3,2}|\leq 6\|F(x)\|_{L^\infty}^4\|x\|_{C^{2,\frac12}}^2(|\dg^3 x|+\int |\eta|^{-\frac12}|\dg^3 x'|d\eta),
$$
proving that $E_{3,2}\in C([0,T];L^2)$. The lower order term $E_{3,3}$ can be estimated similarly and it is also in the same space. Next we continue rewriting
$$
E_{3,4}=\frac{\dg x}{|\dg x|^2}\cdot\int(\dg x_--\dg^2x\eta)\dg^2\Big(\frac{x_-\cdot\dg x_-}{|x_-|^3}\Big)d\eta,
$$
form where we obtain with the same methods the bound
$$
|E_{3,4}|\leq 2\|F(x)\|_{L^\infty}^3\|x\|_{C^{2,\frac12}}(|\dg^3 x|+\int |\eta|^{-\frac12}|\dg^3 x'|d\eta)+C(t).
$$

 It remains to estimate $E_{3,1}$ which can be rewritten as follows
$$
E_{3,1}=\frac{1}{|\dg x|^2}\int\frac{\dg x_-\cdot\dg^4 x'}{|x_-|}d\eta
-\frac{1}{|\dg x|^2}\int\frac{(\dg x\cdot\dg^4 x)_-}{|x_-|}d\eta,
$$
suggesting the splitting
$E_{3,1}=E_{3,1}^1+E_{3,1}^2+E_{3,1}^3+E_{3,1}^4$
where
$$
E_{3,1}^1=\frac{1}{|\dg x|^2}\int\Big(\frac{\dg x_-}{|x_-|}-\frac{\dg^2 x}{|\dg x|}\frac{\eta}{|\eta|}\Big)\cdot\dg^4 x'd\eta,\quad
E_{3,1}^2=\frac{\dg^2x}{|\dg x|^3}\cdot\int\dg^4 x'\frac{\eta}{|\eta|}d\eta,
$$
and
$$
E_{3,1}^3=\frac{-1}{|\dg x|^2}\int (\dg x\cdot\dg^4 x)_- \Big(\frac{1}{|x_-|}-\frac{1}{|\dg x||\eta|}\Big)d\eta,\quad
E_{3,1}^4=\frac{-1}{|\dg x|^3}\int \frac{ (\dg x\cdot\dg^4 x)_-}{|\eta|} d\eta.
$$
We have the kernel:
$$
\Omega_4=\frac{\dg x_-}{|x_-|}-\frac{\dg^2 x}{|\dg x|}\frac{\eta}{|\eta|}=\frac{\dg x_--\dg^2 x\eta}{|x_-|}+
\dg^2 x\eta\Big(\frac{1}{|x_-|}-\frac{1}{|\dg x||\eta|}\Big)
$$
and
$$
E_{3,1}^1=\frac{-1}{|\dg x|^2}\int\de\Omega_4\cdot\dg^3 x'd\eta+\Omega_4\dg^3 x'\Big|_{\eta=-\pi}^{\eta=\pi}.
$$
Dealing with $\de\Omega_4$ in a similar manner as we did before, we get the estimate $|\de\Omega_4|\leq C(t)|\eta|^{-\frac12}$, implying that $E_{3,1}^1\in C([0,T];L^2)$. 

A convenient integration yields
$$
E_{3,1}^2=\frac{\dg^2x}{|\dg x|^3}\cdot(-2\dg^3 x+\dg^3 x(\g+\pi)+\dg^3 x(\g-\pi)),
$$ 
from where the appropriate estimate for $E_{3,1}^2$ follows. Identity \eqref{unacancelacion} allows to obtain
$$
E_{3,1}^3=\frac{3}{|\dg x|^2}\int (\dg^2 x\cdot\dg^3 x)_- \Big(\frac{1}{|x_-|}-\frac{1}{|\dg x||\eta|}\Big)d\eta,
$$
and therefore
$$
|E_{3,1}^3|\leq C\|F(x)\|_{L^\infty}^4\|\dg^2x\|^2_{L^\infty}\int (|\dg^3 x|+|\dg^3 x'|)d\eta.
$$
Finally, using one more time \eqref{unacancelacion} we get
$$
E_{3,1}^4=\frac{3}{|\dg x|^3}\int \frac{ (\dg^2 x\cdot\dg^3 x)_-}{|\eta|} d\eta=\frac{-2}{3}E_{2,1}^3,
$$
where $E_{2,1}^3$ is given in \eqref{e213}. Then, $E_{3,1}^4$ can also be estimated as before. We are done with $E_{3,4}$ and therefore with $E_{3}$. It gives $\lambda(x)\in C([0,T];H^{\frac{3}{\log}})$ as desired. 

Regarding $K$ in \eqref{JK}, we have $K=K_1+K_2+K_3$ where
$$
K_1=\int\dl\dg^2\phi\dl(\dg^2\lambda(x)\dg\phi)d\g,\quad K_2=\int\dl\dg^2\phi\dl(\dg\lambda(x)\dg^2\phi)d\g,
$$
and
$$
K_3=\int\dl\dg^2\phi\dl\dg(\lambda(x)\dg^2\phi)d\g.
$$
At this point it is easy to get
$$
K_1\leq C\|\dl\dg^2\phi\|_{L^2}\|\dl\dg^2\lambda(x)\|_{L^2}\|\dl\dg\phi\|_{L^2}\leq C(t)\|\dl\dg^2\phi\|^2_{L^2},
$$
and
$$
K_2\leq  C\|\dl\dg^2\phi\|^2_{L^2}\|\dl\dg\lambda(x)\|_{L^2}\leq C(t)\|\dl\dg^2\phi\|^2_{L^2}.
$$
For $K_3$ the commutator estimate \eqref{conmutador} allow us to get
\begin{align*}
K_3
&=\int\dl\dg^2\phi(\dl\dg(\lambda(x)\dg^2\phi)-\lambda(x)\dl\dg\dg^2\phi)d\g+
\int\lambda(x)\dl\dg^2\phi\dg\dl\dg^2\phi d\g\\
&\leq C\|\dl\dg^2\phi\|_{L^2}(\|\lambda(x)\|_{H^2}\|\dl\dg^2\phi\|_{L^2}+\|\dl\dg\lambda(x)\|_{L^2}
\|\dg^2\phi\|_{H^{\frac56}})-\frac12\int\dg\lambda(x)|\dl\dg^2\phi|^2 d\g
\end{align*}
to obtain finally
$$
K_3\leq  C\|\lambda(x)\|_{H^2}\|\dl\dg^2\phi\|^2_{L^2}\leq C(t) \|\dl\dg^2\phi\|^2_{L^2}.
$$
Having such good estimates for $K$ and $J$ we can go back to \eqref{JK} and obtain
$$
\frac{d}{dt}\|\dl\dg^2\phi\|_{L^2}^2\leq C(t)\|\dl\dg^2\phi\|^2_{L^2},
$$
which together with \eqref{phiL2} yields
$$
\frac{d}{dt}\|\phi-id\|_{H^{\frac{3}{\log}}}^2\leq C(t)\|\phi-id\|^2_{H^{\frac{3}{\log}}}+C(t),
$$
and then Gronwall Lemma gives existence so long as $\D\int_0^tC(s)ds<\infty$.

Uniqueness then follows similarly because we have
$$
\frac{d}{dt}\|\phi^2-\phi^1\|_{L^2}^2\leq C(t)\|\phi^2-\phi^1\|^2_{L^2}
$$
where $\phi^2$ and $\phi^1$ are two solutions of the equation and $\phi^2(x,0)=\phi^1(x,0)$, and because the above inequality can be obtained with the method described before.

It remains to show that $\dg\phi(\g,t)>0$ for some positive time. This is done with the observation
$$
\dg\phi(\g,t)=\dg\phi(\g,0)+\int_0^t \dg\phi_t(\g,s)ds\geq \min_{\g}\dg\phi(\g,0)-\int_0^t |\dg\phi_t(\g,s)|ds.
$$
The fact that $|\dg\phi_t(\g,s)|\leq C(t)\|\phi\|_{H^{\frac{3}{\log}}}$ implies that $\phi$ remains as a legitimate change of variable so long as
$$
\min_{\g}\dg\phi(\g,0)>\int_0^t  C(s)\|\phi\|_{H^{\frac{3}{\log}}}(s)ds.
$$
%%%%%%%%%%%%%%%%%%%%%%%%%%%%%%%%%%%%%%%%%%%%%%%%%%%%%%%%%%%%%%%%%%%%%%%%%%%%%%%%%%%%%%%%%%%%%%%%%%%%%%%%%%%%%%%%%%%%%%
%%%%%%%%%%%%%%%%%%%%%%%%%%%%%%%%%%%%%%%%%%%%%%%%%%%%%%%%%%%%%%%%%%%%%%%%%%%%%%%%%%%%%%%%%%%%%%%%%%%%%%%%%%%%%%%%%%%%%%
%%%%%%%%%%%%%%%%%%%%%%%%%%%%%%%%%%%%%%%%%%%%%%%%%%%%%%%%%%%%%%%%%%%%%%%%%%%%%%%%%%%%%%%%%%%%%%%%%%%%%%%%%%%%%%%%%%%%%%

\section{Uniqueness for the SQG patch problem}

This section is devoted to show the proof of uniqueness of SQG weak solutions given by patches. In order to do that we introduce the following notation for simply connected domains. 

We say that a bounded simply connected domain $D\subset \R^2$ is $C^{2,\delta}(\T)$ for $0<\delta<1$ if there exists a parameterization of the boundary
$$
\partial D=\{x(\g)\in\R^2: \g\in\T,\, 2\pi\mbox{-periodic}\}
$$
such that $x(\g)\in C^{2,\delta}(\T)$. In particular, a domain $\Omega\in C^{2,\delta}(\T)$ given by
$$
\partial \Omega=\{y(\xi)\in\R^2: \xi\in\T,\, 2\pi\mbox{-periodic}\}
$$
is said to be equal to $D$ if there exists a change of variable
$$
\varphi:\T\to\T,\quad \mbox{biyective},\quad \varphi'(\gamma)>0,\quad \varphi(\g)-\g \quad 2\pi\mbox{-periodic},\quad \varphi\in C^{2,\delta}(\T),
$$ 
such that $x(\g)=y(\varphi(\g))$. Furthermore, a time dependent simply connected domain $D(t)$ belongs to $C([0,T];C^{2,\delta}(\T))\cap C^1([0,T];C^1(\T))$ if there exist parameterizations of the boundaries 
$$
\partial D(t)=\{x(\g,t)\in\R^2: \g\in\T,\, t\in[0,T],\, 2\pi\mbox{-periodic in }\gamma\}
$$	
such that $x(\g,t)\in C([0,T];C^{2,\delta}(\T))\cap C^1([0,T];C^1(\T))$.

The main result in the section is the following.

\begin{thm}\label{USQGp}
	  Consider a solution of \eqref{weaksolution} with $\theta(x,t)$ given by a patch \eqref{patchsolution} and $D^{j}(t)$ a time dependent simply connected domain in $C([0,T];C^{2,\delta}(\T))\cap C^1([0,T];C^1(\T))$ whose free boundary satisfies the arc-chord condition for any $t\in [0,T]$. Furthermore, assume that the function $\bar{\theta}(x,t)$ given by
	  	\begin{equation*}
	  	\bar{\theta}(x,t)=\left\{\begin{array}{rl}
	  	\theta^1,& x\in \bar{D}^1(t),\\
	  	\theta^2,& x\in \bar{D}^2(t)=\R^2\setminus \bar{D}^1(t),
	  	\end{array}\right. 
	  	\end{equation*}
	  	satisfies \eqref{weaksolution} with $\bar{D}^{j}(t)\in C([0,T];C^{2,\delta}(\T))\cap C^1([0,T];C^1(\T))$ and $\theta(x,0)=\bar{\theta}(x,0)$. Then $\theta(x,t)=\bar{\theta}(x,t)$ for any $t\in [0,T]$.
\end{thm}

\emph{Proof:} We consider a solution $\theta(x,t)$ satisfying the hypothesis above. Then, it is shown in \cite{FG}, the parameterization of the free boundary has to fulfill equation \eqref{sqgpatchnormal} where, without loss of generality, we can assume that $\theta_2-\theta_1=\pi$. The length of the curve is
$$
l(t)=\int_{\T} |\dg x(\g,t)|d\g,
$$  
and we shall consider the following change of variable 
$$
\phi(\cdot,t):\T\to\T,\quad \phi(\g,t)=-\pi+\frac{2\pi}{l(t)}\int_{-\pi}^{\g}|\dg x(\eta,t)|d\eta.
$$
Consequently we get the reparameterization
$$
\tilde{x}(\xi,t)=x(\phi^{-1}(\xi,t),t),\qquad\qquad x(\g,t)=\tilde{x}(\phi(\g,t),t),\qquad \xi=\phi(\g,t),
$$
satisfying property \eqref{unacancelacion} ( 
$
|\partial_{\xi}\tilde{x}(\xi,t)|=(2\pi)^{-1}l(t)
$) and having the same regularity ($\tilde{x}(\xi,t)\in C([0,T],C^{2,\delta})\cap C^1([0,T];C^1(\T)))$. As we point out before, the curve $\tilde{x}(\xi,t)$ is a solution of \eqref{sqgpatchnormal} with the tilde notation. We mean by this that $\tilde{x}(\xi,t)$ is a solution of \eqref{sqgpatchnormal} replacing $x$ by $\tilde{x}$ and $\gamma$ by $\xi$. 

For this new evolving curve $\tilde{x}$, the  identity 
$$
\tilde{x}_t(\xi,t)=\tilde{x}_t(\xi,t)\cdot\frac{\partial_{\xi}\tilde{x}(\xi,t)}{|\partial_{\xi}\tilde{x}(\xi,t)|}\frac{\partial_{\xi}\tilde{x}(\xi,t)}{|\partial_{\xi}\tilde{x}(\xi,t)|}+\tilde{x}_t(\xi,t)\cdot\frac{\partial^\bot_{\xi}\tilde{x}(\xi,t)}{|\partial_{\xi}\tilde{x}(\xi,t)|}\frac{\partial^\bot_{\xi}\tilde{x}(\xi,t)}{|\partial_{\xi}\tilde{x}(\xi,t)|}
$$
together with \eqref{sqgpatchnormal} provides
$$
\tilde{x}_t(\xi,t)=\tilde{\mu}(\xi,t)\frac{\partial_{\xi}\tilde{x}(\xi,t)}{|\partial_{\xi}\tilde{x}(\xi,t)|}+\int_{\T}\frac{(\partial_{\xi}\tilde{x}(\xi,t)-\partial_{\xi} \tilde{x}(\xi-\zeta,t))}{|\tilde{x}(\xi,t)-\tilde{x}(\xi-\zeta,t)|}d\zeta \cdot\frac{\partial_{\xi}^{\bot}x(\xi,t)}{|\partial_{\xi}\tilde{x}(\xi,t)|}
\frac{\partial^\bot_{\xi}\tilde{x}(\xi,t)}{|\partial_{\xi}
\tilde{x}(\xi,t)|},
$$
where we have defined $\tilde{\mu}(\xi,t)=\tilde{x}_t(\xi,t)\cdot\partial_{\xi}\tilde{x}(\xi,t)/|\partial_{\xi}\tilde{x}(\xi,t)|$. Taking 
$$
\tilde{\mu}(\xi,t)=\int_{\T}\frac{(\partial_{\xi}\tilde{x}(\xi,t)-\partial_{\xi} \tilde{x}(\xi-\zeta,t))}{|\tilde{x}(\xi,t)-\tilde{x}(\xi-\zeta,t)|}d\zeta \cdot\frac{\partial_{\xi}\tilde{x}(\xi,t)}{|\partial_{\xi}\tilde{x}(\xi,t)|}+\tilde{\lambda}(\tilde{x})(\xi,t)|\partial_{\xi}\tilde{x}(\xi,t)|,
$$
it is easy to find that $\tilde{x}$ satisfies \eqref{QGm} with the tilde notation:
$$
\tilde{x}_t(\xi,t)=\int_{\T}\frac{\partial_{\xi}\tilde{x}(\xi,t)-\partial_{\xi} \tilde{x}(\xi-\zeta,t)}{|\tilde{x}(\xi,t)-\tilde{x}(\xi-\zeta,t)|}d\zeta +\tilde{\lambda}(\tilde{x})(\xi,t)\partial_{\xi}\tilde{x}(\xi,t),
$$
where
$$
\tilde{\lambda}(\tilde{x})(\xi,t)=\Big(\tilde{x}_t(\xi,t)-\int_{\T}\frac{\partial_{\xi}\tilde{x}(\xi,t)-\partial_{\xi} \tilde{x}(\xi-\zeta,t)}{|\tilde{x}(\xi,t)-\tilde{x}(\xi-\zeta,t)|}d\zeta\Big)\cdot\frac{\partial_{\xi}\tilde{x}(\xi,t)}{|\partial_{\xi}\tilde{x}(\xi,t)|^2}.
$$

The regularity of $\tilde{x}(\xi,t)$ yields $\tilde{\lambda}(\tilde{x})(\xi,t)\in C([0,T];C^1(\T))$. Then we can find a function $a\in C^1([0,\tilde{T}];\R)$ as a unique solution of the o.d.e.
$$
a'(t)=\tilde{\lambda}(\tilde{x})(-\pi-a(t),t),\quad a(0)=0,
$$
where $0<\tilde{T}$ by the Picard-Lindel\"of theorem. Since $\sup_{[0,T]}\|\tilde{\lambda}(\tilde{x})\|_{L^\infty}(t)\leq C_m(x),$
for $C_m(x)$ depending on $\sup_{[0,T]}\|F(x)\|_{L^\infty}(t)$ and $\sup_{[0,T]}(\|x\|_{C^{2,\delta}}(t)+\|x_t\|_{L^\infty}(t))$, the function $a(t)$ can be extended to $[0,T]$ satisfying that $|a(t)|\leq TC_m(x)$ for any $t\in [0,T]$.

 The new curve given by $\tilde{x}(\xi,t)=\bar{x}(\xi+a(t),t)$ satisfies 
$$
\bar{x}_t(\al,t)=\int_{\T}\frac{\partial_{\al}\bar{x}(\al,t)-\partial_{\al}\bar{x}(\al\!-\!\beta,t)}{|\bar{x}(\al,t)-\bar{x}(\al\!-\!\beta,t)|}d\beta +\bar{\lambda}(\bar{x})(\al,t)\partial_{\al}\bar{x}(\al,t),
$$
for $\al=\xi-a(t)$ and $\bar{\lambda}(\bar{x})(\al,t)=\tilde{\lambda}(\tilde{x})(\al\!-\!a(t),t)-
\tilde{\lambda}(\tilde{x})(-\!\pi\!-\!a(t),t)$. Since $\partial_{\al} |\partial_{\al} x(\al,t)|=0$ and $\bar{\lambda}(\bar{x})(-\pi,t)=0$, we proceed as in \cite{FG} (see pg. 2585) to find that $\bar{x}$ evolves according to equations (\ref{QGm},\ref{la}) replacing $x$ by $\bar{x}$ and $\gamma$ by $\al$. In particular it is easy to check that $\bar{x}(\al,t)$ has the same regularity than $\tilde{x}(\xi,t)$ and $\tilde{x}(\xi,0)=\bar{x}(\xi,0)$.

 We consider next another solution $\bar{\theta}(x,t)$, satisfying the hypothesis above with the free boundary parameterized by $y(\g,t)\in C([0,T];C^{2,\delta}(\T))\cap C^1([0,T];C^1(\T))$. As $\theta(x,0)=\bar{\theta}(x,0)$, we use a function $\varphi\in C^{2,\delta}(\T)$ to define $\breve{y}(\g,t)=y(\varphi(\g),t)$ in such a way that $\breve{y}(\g,0)=x(\g,0)$. Therefore, it is easy to see that $\breve{y}$ has the same regularity than $y$ and fulfills equation \eqref{sqgpatchnormal}, providing the free boundary of the same patch solution $\bar{\theta}(x,t)$. Next, we reparameterize $\breve{y}(\g,t)$ as we did for $x(\g,t)$ to get $\tilde{y}(\xi,t)$ satisfying  $\partial_{\xi}(|\partial_\xi \tilde{y}(\xi,t)|)=0$ and $\tilde{y}(\xi,0)=\tilde{x}(\xi,0)$. Then we obtain $\bar{y}(\alpha,t)$ similarly as before providing us a solution of equations (\ref{QGm},\ref{la}) after replacing $x$ by $\bar{y}$ and $\gamma$ by $\al$. In particular, all this reparameterization process provides $
 \bar{y}(\alpha,t)$ with the same kind of regularity and satisfying $\bar{x}(\alpha,0)=\bar{y}(\alpha,0)$.

From now on, we will drop the bars for simplicity, using the variables $\g$ and $\eta$ instead of $\alpha$ and $\beta$. As before we shall write $f=f(\g,t)$, $f'=f(\g-\e,t)$, $f_-=f-f'$ and $\int=\int_{\T}$, when there is no danger of confusion in the writing of our double integrals in variables $\g$ and $\eta$. During the time of existence $T>0$ one has the arc-chord condition $F(x)$  in $L^\infty(0,T;L^{\infty})$. In the following $C$ will denote a constant which may be different from inequality to inequality but  depending only on $\sup_{[0,T]}\|x\|_{C^{2,\delta}}(t)$, $\sup_{[0,T]}\|y\|_{C^{2,\delta}}(t)$, $\sup_{[0,T]}\|F(x)\|_{L^{\infty}}(t)$ and $\sup_{[0,T]}\|F(y)\|_{L^{\infty}}(t)$.

 Let us consider the function $z(\g,t)=x(\g,t)-y(\g,t)$, we have
$$
\frac{1}{2}\frac{d}{dt}\|z\|^2_{L^2}=\int z\cdot z_t d\g= I_1+I_2,
$$
where
$$
I_1=\int z\cdot \int \Big(\frac{\dg x_-}{|x_-|}-\frac{\dg y_-}{|y_-|} \Big)d\e d\g,\quad
I_2=\int z\cdot (\lambda(x)\dg x-\lambda(y)\dg y) d\g.
$$
Let us split $I_1$:
$$
I_{1,1}=\int z\cdot \int \frac{\dg z_-}{|x_-|}d\e d\g,\quad
I_{1,2}=\int z\cdot \int \dg y_-\Big(\frac{1}{|x_-|}-\frac{1}{|y_-|} \Big)d\e d\g.
$$
Then with an adequate change of variables, we obtain
$$
I_{1,1}=\int\int\frac{z(\g)\cdot(\dg z(\g)-\dg z(\e))}{|x(\g)-x(\e)|}d\e d\g=-\int\int\frac{z(\e)\cdot(\dg z(\g)-\dg z(\e))}{|x(\g)-x(\e)|}d\e d\g
$$
thus
$$
I_{1,1}=\frac12\int\int\frac{(z(\g)-z(\e))\cdot(\dg z(\g)-\dg z(\e))}{|x(\g)-x(\e)|}d\e d\g=
\frac12\int\int\frac{z_-\cdot \dg z_-}{|x_-|} d\g d\e.
$$
Integration by parts provides
$$
I_{1,1}=\frac{1}{4}\int\int|z_-|^2\frac{(x_-\cdot\dg x_-)}{|x_-|^3} d\g d\e=
\frac{1}{4}\int\int|z_-|^2F(x)^3\frac1{|\e|} \frac{(x_-\cdot\dg x_-)}{\e^2} d\g d\e
$$
The inequality
\begin{equation}\label{dm}
|(x_-\cdot\dg x_-)-\dg x\cdot\dg^2 x\e^2|\leq 2\|x\|^2_{C^{2,\delta}}|\e|^{2+\delta},
\end{equation}
together with the fact that $\dg x\cdot\dg^2 x=0$ allows us to get
$$
I_{1,1}\leq \|F(x)\|^3_{L^\infty}\|x\|^2_{C^{2,\delta}}\int|\e|^{\delta-1} \int (|z|^2+|z'|^2) d\g d\e \leq C\|z\|^2_{L^2}.
$$
For $I_{1,2}$ one writes
$$
I_{1,2}=\int z\cdot \int \dg y_- \frac{|y_-|-|x_-|}{|x_-||y_-|}d\e d\g\leq \int\int \frac{|z||\dg y_-||z_-|}{|x_-||y_-|}d\e d\g,
$$
which yields
$$
I_{1,2}\leq\int\int |z|\frac{|\dg y_-|}{|\e|}\frac{|z_-|}{|\e|}F(x)F(y)d\g d\e.
$$
Then the identity
\begin{equation}\label{tvm} f_-=\e \int_0^1\dg f(\g+(s-1)\e)ds\end{equation}
allows us to get the bound
$$
I_{1,2}\leq \|F(x)\|_{L^\infty}\|F(y)\|_{L^\infty}\|\dg^2 y\|_{L^\infty}\int_0^1\int\int |z| |\dg z(\g+(s-1)\e)| d\g d\e ds,
$$
which yields the desired control: $I_{1,2}\leq C \|z\|^2_{H^1}$. 

Regarding $I_{2}$ we split further
$$
I_{2,1}=\int z\cdot \dg z \lambda(x) d\g,\quad I_{2,2}=\int z\cdot\dg y  (\lambda(x)-\lambda(y)) d\g.
$$
It is easy to get
$$
I_{2,1}\leq \|z\|_{L^2}\|\dg z\|_{L^2}\|\lambda(x)\|_{L^\infty}\leq C \|z\|^2_{H^1},
$$ thus we are done with $I_{2,1}$.

For the reminder term we have
$$
I_{2,2}\leq \|z\|_{L^2}\|\dg y\|_{L^\infty}\|\lambda(x)-\lambda(y)\|_{L^2}\leq C \|z\|_{L^2}\|\lambda(x)-\lambda(y)\|_{L^2},
$$
let us write $\lambda(x)-\lambda(y)=G_1+G_2$ where
$$
G_1=\frac{\g+\pi}{2\pi}\int\Big[\frac{\dg x}{|\dg x|^2}\cdot \dg \Big(\int
\frac{\dg x_-}{|x_-|}d\e \Big)-\frac{\dg y}{|\dg y|^2}\cdot \dg \Big(\int
\frac{\dg y_-}{|y_-|}d\e \Big)\Big] d\g,
$$
and
\begin{align*}
G_2=&-\int_{-\pi}^\g \frac{\dg x(\e,t)}{|\dg x(\e,t)|^2}\cdot \de \Big(\int\frac{\dg
x(\e,t)-\dg x(\e-\xi,t)}{|x(\e,t)-x(\e-\xi,t)|}d\xi \Big)d\e\\
&+\int_{-\pi}^\g \frac{\dg y(\e,t)}{|\dg y(\e,t)|^2}\cdot \de \Big(\int\frac{\dg
y(\e,t)-\dg y(\e-\xi,t)}{|y(\e,t)-y(\e-\xi,t)|}d\xi \Big)d\e.
\end{align*}
Then we decompose further $G_1=G_{1,1}+G_{1,2}+G_{1,3}+G_{1,4}$:
$$
G_{1,1}=\frac{\g+\pi}{2\pi}\int \frac{\dg z}{|\dg x|^2}\cdot \dg \Big(\int
\frac{\dg x_-}{|x_-|}d\e \Big)d\g,
$$
$$
G_{1,2}=\frac{\g+\pi}{2\pi}\int \big(\frac{1}{|\dg x|^2}-\frac{1}{|\dg y|^2}\big)\dg y\cdot \dg \Big(\int
\frac{\dg x_-}{|x_-|}d\e \Big)d\g,
$$
$$
G_{1,3}=-\frac{\g+\pi}{2\pi}\int \frac{\dg y}{|\dg y|^2}\cdot \int
\Big(\dg x_- \frac{x_-\cdot \dg x_-}{|x_-|^3}-\dg y_-\frac{y_-\cdot \dg y_-}{|y_-|^3}\Big)d\e d\g,
$$
and
$$
G_{1,4}=\frac{\g+\pi}{2\pi}\int \frac{\dg y}{|\dg y|^2}\cdot \int
\Big(\frac{\dg^2 x_-}{|x_-|}-\frac{\dg^2 y_-}{|y_-|}\Big)d\e d\g.
$$
We proceed as before $$|G_{1,1}|\leq C(\|F(x)\|^3_{L^\infty}\|\dg^2 x\|_{C^\delta}+\|F(x)\|^4_{L^\infty}\|\dg^2 x\|^2_{L^\infty})\|\dg z\|_{L^2}$$
and therefore $\|G_{1,1}\|_{L^2}\leq C\|z\|_{H^1}$. In a similar way we find $\|G_{1,2}\|_{L^2}\leq C\|z\|_{H^1}$. To estimate $G_{1,3}$ we write $G_{1,3}=G_{1,3,1}+G_{1,3,2}+G_{1,3,3}$ where $G_{1,3,1}$ and $G_{1,3,2}$ are the most singular terms:
$$
G_{1,3,1}=-\frac{\g+\pi}{2\pi}\int \frac{\dg y}{|\dg y|^2}\cdot \int
\dg z_-\frac{x_-\cdot \dg x_-}{|x_-|^3}d\e d\g,
$$
$$
G_{1,3,2}=-\frac{\g+\pi}{2\pi}\int \frac{\dg y}{|\dg y|^2}\cdot \int
\dg y_-\frac{x_-\cdot \dg z_-}{|x_-|^3}d\e d\g,
$$
because $G_{1,3,3}$ satisfies obviously the desired bound: $\|G_{1,3,3}\|_{L^2}\leq C\|z\|_{H^1}$. To control $I_{1,1}$, we use \eqref{dm} and the fact that $\dg x\cdot\dg^2 x=0$ that is:
$$
|G_{1,3,1}|\leq \|F(y)\|_{L^\infty}\|F(x)\|^3_{L^\infty}\|x\|^2_{C^{2,\delta}}\int |\e|^{\delta-1}\int (|\dg z|+|\dg z'|)d\g d\e,
$$
implying $\|G_{1,3,1}\|_{L^2}\leq C\|z\|_{H^1}$. 

Inside the expression of $G_{1,3,2}$ we observe that
$$
\dg y\cdot \dg y_-=\dg y\cdot (\dg y_--\eta\dg^2y)
$$
which together with the estimate
\begin{equation}\label{h}
|\dg y_--\eta\dg^2y|\leq \|y\|_{C^{2,\delta}}|\eta|^{1+\delta},
\end{equation}
give us
$$
|G_{1,3,2}|\leq \|F(y)\|_{L^\infty}\|F(x)\|^2_{L^\infty}\|y\|^2_{C^{2,\delta}}\int |\e|^{\delta-1}\int (|\dg z|+|\dg z'|)d\g d\e,
$$
and $\|G_{1,3,2}\|_{L^2}\leq C\|z\|_{H^1}$. 

Next let us write $G_{1,4}=G_{1,4,1}+G_{1,4,2}$ where
$$
G_{1,4,1}=\frac{\g+\pi}{2\pi}\int \frac{\dg y}{|\dg y|^2}\cdot \int \dg^2 x_-
\Big(\frac{1}{|x_-|}-\frac{1}{|y_-|}\Big)d\e d\g,
$$
$$
G_{1,4,2}=\frac{\g+\pi}{2\pi}\int \frac{\dg y}{|\dg y|^2}\cdot \int
\frac{\dg^2 z_-}{|y_-|}d\e d\g.
$$
Equality \eqref{tvm} allows us to obtain
$$
|G_{1,4,1}|\leq \|F(y)\|^2_{L^\infty}\|F(x)\|_{L^\infty}\|\dg^2 x\|_{C^\delta}\int_0^1\int|\eta|^{\delta-1}\int|\dg z(\g+(r-1)\e)|d\g d\e dr,
$$
and hence $\|G_{1,4,1}\|_{L^2}\leq C\|z\|_{H^1}$. Integration by parts allows us to decompose further $G_{1,4,2}=G_{1,4,2}^1+G_{1,4,2}^2$ where
$$
G^1_{1,4,2}=\frac{\g+\pi}{2\pi}\int\! \frac{\dg y}{|\dg y|^2}\cdot \int
\frac{\dg z_- (y_-\cdot \dg y_-)}{|y_-|^{3}}d\e d\g,
\quad
G^2_{1,4,2}=-\frac{\g+\pi}{2\pi}\int\! \frac{\dg^2 y}{|\dg y|^2}\cdot \int
\frac{\dg z_-}{|y_-|}d\e d\g.
$$
The first term can be estimated as $G_{1,3,1}$:
$$
|G^1_{1,4,2}|\leq C\|F(y)\|^4_{L^\infty}\|y\|^2_{C^{2,\delta}}\|z\|_{H^1}\leq C\|z\|_{H^1}.
$$
We symmetrize $G^2_{1,4,2}$ as in $I_{1,1}$:
$$
G^2_{1,4,2}=-\frac{\g+\pi}{4\pi|\dg y|^2}\int\int
\frac{\dg^2 y_-\cdot\dg z_-}{|y_-|}d\e d\g
$$
which yields the estimate: 
$$|G^2_{1,4,2}|\leq C \|F(y)\|^3_{L^\infty}\|\dg^2 y\|_{C^\delta}\|\dg z\|_{L^2}\leq C\|z\|_{H^1},
$$ implying that
$$|G_1|\leq  C\|z\|_{H^1}.$$

For the sake of simplicity we exchange the variables in $G_2$ so that
\begin{align*}
G_2=&-\int_{-\pi}^\xi \frac{\dg x}{|\dg x|^2}\cdot \dg \Big(\int\frac{\dg
x_-}{|x_-|}d\e \Big)d\g+\int_{-\pi}^\xi \frac{\dg y}{|\dg y|^2}\cdot \dg \Big(\int\frac{\dg
y_-}{|y_-|}d\e \Big)d\g.
\end{align*}
We claim that $\|G_2\|_{L^2}\leq C\|z\|_{H^1}$. To show that we decompose further $G_{2}=G_{2,1}+G_{2,2}$ where
\begin{align*}
G_{2,1}=&\int_{-\pi}^\xi \frac{\dg x}{|\dg x|^2}\cdot \int\dg
x_-\frac{x_-\cdot\dg
x_-}{|x_-|^3}d\e d\g-\int_{-\pi}^\xi \frac{\dg y}{|\dg y|^2}\cdot \int\dg
y_-\frac{y_-\cdot\dg
y_-}{|y_-|^3}d\e d\g,
\end{align*}
and
\begin{align*}
G_{2,2}=&-\int_{-\pi}^\xi \frac{\dg x}{|\dg x|^2}\cdot\int\frac{\dg^2
x_-}{|x_-|}d\e d\g+\int_{-\pi}^\xi \frac{\dg y}{|\dg y|^2}\cdot \int\frac{\dg^2
y_-}{|y_-|}d\e d\g.
\end{align*}
We deal with $G_{2,1}$ as with $G_1$, to obtain $|G_{2,1}|\leq C\|z\|_{H^1}$. The identities
$$
\dg x\cdot\dg^2x_-=-\dg x\cdot\dg^2x'=-\dg x_-\cdot\dg^2x'
$$
allow us obtain
\begin{align*}
G_{2,2}=&\int_{-\pi}^\xi\int\frac{\dg x_-\cdot\dg^2
x'}{|\dg x|^2|x_-|}d\e d\g-\int_{-\pi}^\xi\int
\frac{\dg y_-\cdot\dg^2 y'}{|\dg y|^2|y_-|}
d\e d\g.
\end{align*}
A new decomposition yields $G_{2,2}=G_{2,2,1}+G_{2,2,2}+G_{2,2,3}$ where
$$
G_{2,2,1}=\int_{-\pi}^\xi\int\frac{\dg z_-\cdot\dg^2
x'}{|\dg x|^2|x_-|}d\e d\g,\quad G_{2,2,2}=\int_{-\pi}^\xi\int\frac{\dg y_-\cdot\dg^2
z'}{|\dg x|^2|x_-|}d\e d\g,
$$
and  $G_{2,2,3}$ collects  the lower order characters, which can be estimate as before: $|G_{2,2,3}|\leq C\|z\|_{H^1}$. One has
$$
G_{2,2,1}=\int_{-\pi}^\xi\Big[\dg\Big(\int\frac{z_-\cdot\dg^2
x'}{|\dg x|^2|x_-|}d\e\Big)-\int z_-\cdot\dg\Big(\frac{\dg^2
x'}{|\dg x|^2|x_-|}\Big)d\e\Big] d\g,
$$
which helps to decompose as follows: $G_{2,2,1}=G_{2,2,1}^1+G_{2,2,1}^2+G_{2,2,1}^3$ where
$$
G_{2,2,1}^1=\Big(\int\frac{z_-\cdot\dg^2
x'}{|\dg x|^2|x_-|}d\e\Big)\Big|_{\gamma=-\pi}^{\gamma=\xi},\quad G_{2,2,1}^2=\int_{-\pi}^\xi\int z_-\cdot\frac{\dg^3
x'}{|\dg x|^2|x_-|}d\e d\g,
$$
and $G_{2,2,1}^3$ consists of the lower order terms. At this point it is easy to get the estimate $|G_{2,2,1}^3|\leq C\|z\|_{H^1}$   
and
\begin{align*}
|G_{2,2,1}^1|\leq&2\|F(x)\|^3_{L^\infty}\|z\|_{C^{\frac12}}\|\dg^2
x\|_{L^\infty}\int|\e|^{-\frac12} d\e\leq  C\|z\|_{H^1}
\end{align*}
as a consequence of Sobolev's embedding. Concerning $G_{2,2,1}^2$ we write

$
\dg^3 x'=\de\dg^2 x_-
$ and integrate by parts to find
$$
G_{2,2,1}^2=\int_{-\pi}^\xi\int z_-\cdot\dg^2 x_-\frac{ x_-\cdot\dg x' }{|\dg x|^2|x_-|^3}d\e d\g
-\int_{-\pi}^\xi\int \dg z'\cdot\frac{\dg^2 x_-}{|\dg x|^2|x_-|}d\e d\g.
$$
Proceeding as before we obtain
\begin{align*}
|G_{2,2,1}^2|\leq& \|F(x)\|^3_{L^\infty}\|\dg^2x\|_{C^\delta}\int_0^1\int|\eta|^{\delta-1}\int
|\dg z(\g+(r-1)\eta)|d\g d\eta dr\\
&+\|F(x)\|^3_{L^\infty}\|\dg^2x\|_{C^\delta}\int|\eta|^{\delta-1}
\int|\dg z'|d\g d\eta \leq C\|z\|_{H^1}.
\end{align*}
Gathering together the last three estimates we have $|G_{2,2,1}|\leq C\|z\|_{H^1}$.

Regarding $G_{2,2,2}$  identity
$
\dg^2 z'=\de\dg z_-
$
and integration by parts yield
$$
G_{2,2,2}=-\int_{-\pi}^\xi\int\frac{\dg^2 y'\cdot\dg
z_-}{|\dg x|^2|x_-|}d\e d\g+\int_{-\pi}^\xi\int\frac{(\dg y_-\cdot\dg
z_-)(x_-\cdot\dg x')}{|\dg x|^2|x_-|^3}d\e d\g.
$$
In the formula above we find two terms analogous to those of $G_{2,2,1}$, so that a similar 
argument gives us $|G_{2,2,2}|\leq C\|z\|_{H^1}$. Thereby we have finally obtained $I_{2,2}\leq C\|z\|^2_{H^1}$.

A consequence of all those estimates is the differential inequalities:
$$
\frac{d}{dt}\|z\|^2_{L^2}\leq C\|z\|^2_{H^1}.
$$

The next step is to analysed
$$
\frac{1}{2}\frac{d}{dt}\|\dg z\|^2_{L^2}=\int \dg z\cdot \dg z_t d\g=I_3+I_4,
$$
where
$$
I_3=\int \dg z\cdot \int \dg\Big(\frac{\dg x_-}{|x_-|}-\frac{\dg y_-}{|y_-|} \Big)d\e d\g,\quad
I_4=\int \dg z\cdot \dg(\lambda(x)\dg x-\lambda(y)\dg y) d\g.
$$
We split further $I_3:$
$$
I_{3,1}=\int \dg z\cdot \int \Big(\frac{\dg^2 x_-}{|x_-|}-\frac{\dg^2 y_-}{|y_-|} \Big)d\e d\g,
$$
$$
I_{3,2}=\int \dg z\cdot \int \Big(-\frac{\dg x_-(x_-\cdot \dg x_-)}{|x_-|^3}+\frac{\dg y_-(y_-\cdot \dg y_-)}{|y_-|^3} \Big)d\e d\g.
$$
Then we write $I_{3,1}=I_{3,1,1}+I_{3,1,2}$ where
$$I_{3,1,1}=\int \dg z\cdot \int \frac{\dg^2 z_-}{|x_-|} d\e d\g,
\quad
I_{3,1,2}=\int \dg z\cdot \int \dg^2 y_-\frac{|y_-|-|x_-|}{|x_-||y_-|} d\e d\g
$$
Replacing in $I_{1,1}$ $z$ by $\dg z$ we find $I_{3,1,1}$, and
$$
I_{3,1,1}=\frac14\int\int |\dg z_-|^2  \frac{x_-\cdot\dg x_-}{|x_-|^3} d\e d\g\leq C\|\dg z\|^2_{L^2}.
$$
At this stage of the proof we can easily obtain the estimate
$$
I_{3,1,2}\leq \|F(x)\|_{L^\infty}\|F(y)\|_{L^\infty}\|\dg^2y\|_{C^{\delta}}\!\int_0^1\!\int\! |\e|^{\delta-1}\!\!\int\! |\dg z||\dg z(\g\!+\!(s\!-\!1)\eta)|d\g d\e ds\leq C\|z\|_{H^1}^2,
$$
and we are done with $I_{3,1}$. 

For $I_{3,2}$ we split further: $I_{3,2}=I_{3,2,1}+I_{3,2,2}+I_{3,2,3}+I_{3,2,4}$ where
$$
I_{3,2,1}=-\int \dg z\cdot \int \frac{\dg z_-(x_-\cdot \dg x_-)}{|x_-|^3} d\e d\g,\quad I_{3,2,2}=-\int \dg z\cdot \int \frac{\dg y_-(z_-\cdot \dg x_-)}{|x_-|^3} d\e d\g,
$$
$$
I_{3,2,3}=-\int \dg z\cdot \int \dg y_-(y_-\cdot \dg x_-)(|x_-|^{-3}-|y_-|^{-3}) d\e d\g
$$
and
$$
I_{3,2,4}=-\int \dg z\cdot \int \frac{\dg y_-(y_-\cdot \dg z_-)}{|y_-|^3} d\e d\g.
$$
Inequality \eqref{dm} yields $I_{3,2,1}\leq C\|z\|^2_{H^1}$. No cancellation is needed to get
$$
I_{3,2,2}\leq C\|z\|^2_{H^1},\quad I_{3,2,3}\leq C\|z\|^2_{H^1}
$$
On the other hand, we pay special attention to $I_{3,2,4}$. By identity \eqref{tvm} we split it further
$$
I_{3,2,4}^1=-\int\int \dg z\cdot\frac{\dg y_-}{|y_-|^3}\int_0^1(y_--\dg y(\g+(r-1)\e)\e)\cdot \dg^2 z(\g+(r-1)\e)\e dr d\e d\g,
$$
$$
I_{3,2,4}^2=-\int\int
\dg z\cdot\frac{\dg y_-}{|y_-|^3} \e^2\int_0^1\dg y(\g+(r-1)\e)\cdot \dg^2 z(\g+(r-1)\e)dr d\e d\g.
$$
In $I_{3,2,4}^1$ we have $\dg^2 z(\g+(r-1)\e)\e=\frac{d}{dr}(\dg z(\g+(r-1)\e))$ and integration by parts in $r$ to find
$$
I_{3,2,4}^{1,1}=-\int\int \dg z\cdot\frac{\dg y_-}{|y_-|^3}\int_0^1\dg^2 y(\g+(r-1)\e)\e^2\cdot \dg z(\g+(r-1)\e) dr d\e d\g,
$$
$$
I_{3,2,4}^{1,2}=-\int\!\!\int \dg z\cdot\frac{\dg y_-}{|y_-|^3}(y_--\dg y\e)\cdot \dg z d\e d\g,
\,\,
I_{3,2,4}^{1,3}=\int\!\!\int \dg z\cdot\frac{\dg y_-}{|y_-|^3}(y_--\dg y'\e)\cdot \dg z'd\e d\g.
$$
Proceeding as before, we obtain the estimate $I_{3,2,4}^{1,j}\leq C\|F(y)\|^3_{L^\infty}\|\dg^2 y\|^2_{L^\infty}\|\dg z\|_{L^2}^2\leq C\|z\|_{H^1}^2$, for $1\leq j\leq 3$. To handle $I_{3,2,4}^2$ we observe that 
\begin{equation}\label{lqsc}\dg y\cdot\dg^2 z=\dg y\cdot\dg^2 x=-\dg z\cdot\dg^2 x\end{equation}
to get
$$
I_{3,2,4}^2=\int\int
\dg z\cdot\frac{\dg y_-}{|y_-|^3} \e^2\int_0^1\dg z(\g+(r-1)\e)\cdot \dg^2 x(\g+(r-1)\e)dr d\e d\g.
$$
Finally  we estimate this term  $I_{3,2,4}^2\leq C\|F(y)\|^3_{L^\infty}\|\dg^2 y\|_{L^\infty}\|\dg^2 x\|_{L^\infty}\|\dg z\|_{L^2}^2\leq C\|z\|_{H^1}^2$, which completes the control of $I_3$ . 

 Next we proceed with a last splitting: $I_4=I_{4,1}+I_{4,2}+I_{4,3}+I_{4,4}$ where
$$
I_{4,1}=\int \dg z\cdot \lambda(x)\dg^2 z d\g,\quad
I_{4,2}=\int \dg z\cdot (\lambda(x)-\lambda(y))\dg^2 y d\g,
$$
$$
I_{4,3}=\int |\dg z|^2\dg\lambda(x) d\g,\quad
I_{4,4}=\int \dg z\cdot\dg y (\dg\lambda(x)-\dg\lambda(y)) d\g.
$$
Integration by parts in $I_{4,1}$ yields: $I_{4,1}\leq \frac12\|\dg z\|_{L^2}^2\|\dg\lambda(x)\|_{L^\infty}\leq C\|z\|^2_{H^1}$ (see the last section for more details on the estimates for $\lambda(x)$).

We have $$I_{4,2}\leq C \|\dg^2 y\|_{L^\infty}\|\dg z\|_{L^2}\|\lambda(x)-\lambda(y)\|_{L^2}\leq C\|z\|^2_{H^1}$$ by similar arguments used for $I_{2,2}$. The control of $I_{4,3}$ follows as in $I_{4,1}$. Finally, integration by parts
$$
I_{4,4}=-\int \dg z\cdot\dg^2 y (\lambda(x)-\lambda(y)) d\g-\int \dg^2 z\cdot\dg y (\lambda(x)-\lambda(y)) d\g,
$$
and identity \eqref{lqsc}
allows us to get the estimate:
$$
I_{4,4}\leq (\|\dg^2 y\|_{L^\infty}+\|\dg^2 x\|_{L^\infty})\|\dg z\|_{L^2} \|\lambda(x)-\lambda(y)\|_{L^2}\leq C\|z\|_{H^1}^2.
$$

Therefore we have obtained
$$
\frac{d}{dt}\|z\|_{H^1}\leq C\|z\|_{H^1},
$$
which allows us the use of Gronwall's inequality to get  uniqueness.

Remark: We have proven the equality $\bar{x}(\al,t)=\bar{y}(\al,t)$. Therefore, undoing the re\-pa\-ra\-me\-te\-ri\-za\-tion process, the patch $\theta$ with a moving boundary given by $x(\g,t)$ is the same as the patch $\bar{\theta}$ described by $y(\g,t)$.

\section{Uniqueness for the system \eqref{psqglz}}

This section is devoted to show uniqueness for the system \eqref{psqglz}. The argument shown below is straight, dealing with the system \eqref{psqglz} without any change of parameterization. As before, to simplify notation we shall write $f=f(\g,t)$, $f'=f(\g-\e)$ and $f-f'=f_-$ when there is no danger of confusion.

We consider two solution for the system \eqref{psqglz}:
$$
x_t(\g,t)=\int_{\T}\frac{\dg x(\g,t)-\dg x(\g-\e,t)}{|x(\g,t)-x(\g-\e,t)|}d\e,
$$
given by $x(\gamma,t)$ and $y(\g,t)$ in the space $C([0,T];H^{\frac3\log}(\T))$ with the same initial data. During the time of existence $T>0$ one finds $F(x)$ and $F(y)$ in $C([0,T];L^{\infty}(\T\times\T))$. Here $C$ denotes a constant which may be different from inequality to inequality but only depends on $\sup_{[0,T]}\|x\|_{H^{\frac3\log}}(t)$, $\sup_{[0,T]}\|y\|_{H^{\frac3\log}}(t)$, $\sup_{[0,T]}\|F(x)\|_{L^{\infty}}(t)$ and $\sup_{[0,T]}\|F(y)\|_{L^{\infty}}(t)$.

Let us consider the function $z(\g,t)=x(\g,t)-y(\g,t)$. One finds
$$
\frac{1}{2}\frac{d}{dt}\|z\|^2_{L^2}=\int z\cdot z_t d\g= I_1+I_2,
$$
where
$$
I_1=\int z\cdot \int \frac{\dg z_-}{|x_-|}d\e d\g,\quad
I_2=\int z\cdot\int \dg y_-\Big(\frac{1}{|x_-|}-\frac{1}{|y_-|}\Big) d\e d\g.
$$
Next we symmetrize $I_1$ and integrate by parts to get
$$
I_1=\frac12\int\int z_-\cdot \frac{\dg z_-}{|x_-|}d\e d\g=\frac{1}{4}\int\int|z_-|^2\frac{x_-\cdot\dg x_-}{|x_-|^3} d\g d\e.
$$
We have the splitting: $I_1=I_{1,1}+I_{1,2}$ where
$$
I_{1,1}=\frac{1}{4}\int\int z\cdot z_- \frac{x_-\cdot\dg x_-}{|x_-|^3}  d\e d\g,\quad\mbox{and}\quad
I_{1,2}=-\frac{1}{4}\int\int z'\cdot z_- \frac{x_-\cdot\dg x_-}{|x_-|^3}  d\e d\g.
$$
Then a simple exchange of variables yields $I_{1,1}=I_{1,2}$.
We have:
$$
I_{1,1}=\frac{1}{4}\int\int z\cdot z_- \Big(\frac{x_-\cdot\dg x_-}{|x_-|^3}-\frac{\dg x\cdot \dg^2 x}{|\dg x|^3|\e|}\Big) d\e d\g+ \frac{1}{4}\int\int z\frac{\dg x\cdot \dg^2 x}{|\dg x|^3} \mathcal{L}(z)  d\g,
$$
hence
\begin{equation}\label{I1U}
	|I_1|\leq 2|I_{1,1}|\leq C\|z\|_{L^2}^2+C\|z\|_{L^2}\|\mathcal{L}(z)\|_{L^2}.
\end{equation}
It remains an estimate for $I_2$. We rewrite
$$
I_2=-\int z\cdot\int \dg y_-\frac{(x_-+y_-)\cdot z_-}{|x_-||y_-|(|x_-|+|y_-|)} d\e d\g
$$ 
and decompose $I_2=I_{2,1}+I_{2,2}$ where
$$
I_{2,1}=-\int z\cdot\int \Big(\frac{\dg y_-  (x_-+y_-)\cdot z_-}{|x_-||y_-|(|x_-|+|y_-|)}-\frac{\dg^2 y  (\dg x+\dg y)\cdot z_-}{|\dg x||\dg y|(|\dg x|+|\dg y|)|\e|}\Big) d\e d\g,
$$ 
and
$$
I_{2,2}=-\int z\cdot \dg^2 y \frac{\dg x+\dg y}{|\dg x||\dg y|(|\dg x|+|\dg y|)}\cdot \mathcal{L}(z) d\g.
$$
As before, we control $I_{2,1}$ and $I_{2,2}$ in the following manner:
$$
I_{2,1}\leq C \int |\e|^{-\frac23}\int |z|(|z|+|z'|)d\g d\e \leq C\|z\|_{L^2}^2, \quad I_{2,2}\leq C\|z\|_{L^2}\|\mathcal{L}(z)\|_{L^2}.
$$
Adding both estimates we obtain the bound for $I_2$, which together with \eqref{I1U} yield
\begin{equation}\label{esUF}
	\frac{d}{dt}\|z\|_{L^2}\leq  C(\|z\|_{L^2}+\|\mathcal{L}(z)\|_{L^2}).
\end{equation}
Next we show that
\begin{equation}\label{miU}
	\|\mathcal{L}(f)\|_{L^2}\leq pC\|f\|_{L^2}^{1-\frac{1}{p}}\|\dg f\|_{L^2}^{1/p}.
\end{equation}
We have
$$
\|\mathcal{L}(z)\|^2_{L^2}\leq C\sum_{k\neq 0} \ln^2(2|k|)|\widehat{f}(k)|^2\leq C\big(\sum_{k\neq 0}|\widehat{f}(k)|^2\big)^{1-\frac1p}
\big(\sum_{k\neq 0}\ln^{2p}(2|k|)|\widehat{f}(k)|^2\big)^{\frac{1}{p}},
$$
for $1\leq p<\infty$ and therefore inequality $\ln^p|k|\leq p^p |k|$ with $|k|\geq 1$ gives \eqref{miU}. Introducing that estimate in \eqref{esUF} we obtain 
$$
\frac{d}{dt}\|z\|_{L^2}\leq C p\|z\|^{1-\frac1p}_{L^2},
$$
for $p\geq1$. Since $\|z\|_{L^2}(0)=0$, we can conclude that the maximal solution of this inequality satisfies
$$
\|z\|_{L^2}(t)\leq (C t)^p,
$$
for $p\geq 1$. Therefore, choosing $t\leq (2C)^{-1}$ and taking the limit as $p\rightarrow +\infty$ we prove uniqueness.

\subsection*{{\bf Acknowledgements}}

\smallskip
AC were partially supported by the grant MTM2014-56350-P (Spain). AC and DC were partially supported by the ICMAT Severo Ochoa project SEV-2015-556.
DC and FG were partially supported by the grant MTM2014-59488-P (Spain). 
FG were partially supported by the Ram\'on y Cajal program RyC-2010-07094, the grant P12-FQM-2466 from Junta de Andaluc\'ia (Spain) and the ERC Starting Grant 639227.

%%%%%%%%%%%%%%%%%%%%%%%%%%%%%%%%%%%%%%%%%%%%%%%%%%%%%%%%%%%%%%%%%%%%%%%%%%%%%%%%%%%%%%%%%%%%%%%%%%%%%%%%%%%%%%%%%%%%%
%%%%%%%%%%%%%%%%%%%%%%%%%%%%%%%%%%%%%%%%%%%%%%%%%%%%%%%%%%%%%%%%%%%%%%%%%%%%%%%%%%%%%%%%%%%%%%%%%%%%%%%%%%%%%%%%%%%%%

%%%%%%%%%%%%%%%%%%%%%%%%%%%%%%%%%%%%%%%%%%%%%%%%%%%%%%%%%%%%%%%%%%%%%%%%%%%%%%%%%%%%%%%%%%%%%%%%%%%%%%%%%%%%%%
%%%%%%%%%%%%%%%%%%%%%%%%%%%%%%%%%%%%%%%%%%%%%%%%%%%%%%%%%%%%%%%%%%%%%%%%%%%%%%%%%%%%%%%%%%%%%%%%%%%%%%%%%%%%%%
%%%%%%%%%%%%%%%%%%%%%%%%%%%%%%%%%%%%%%%%%%%%%%%%%%%%%%%%%%%%%%%%%%%%%%%%%%%%%%%%%%%%%%%%%%%%%%%%%%%%%%%%%%%%%%

\begin{quote}
\begin{tabular}{l}
\textbf{Antonio C\'ordoba} \\
{\small Instituto de Ciencias Matem\'aticas (ICMAT)}\\
{\small \& Departamento de Matem\'aticas}\\{\small Facultad de
Ciencias} \\ {\small Universidad Aut\'onoma de Madrid}
\\ {\small Crta. Colmenar Viejo km.~15,  28049 Madrid,
Spain} \\ {\small Email: antonio.cordoba@uam.es}
\end{tabular}
\end{quote}
\begin{quote}
\begin{tabular}{ll}
\textbf{Diego C\'ordoba} &  \textbf{Francisco Gancedo}\\
{\small Instituto de Ciencias Matem\'aticas (ICMAT)} & {\small Departamento de An\'alisis Matem\'atico \& IMUS}\\
{\small Consejo Superior de Investigaciones Cient\'ificas} & {\small Universidad de Sevilla}\\
{\small C/ Nicol\'as Cabrera, 13-15,} & {\small C/ Tarfia, s/n}\\
{\small Campus Cantoblanco UAM, 28049 Madrid, Spain} & {\small Campus Reina Mercedes, 41012, Sevilla, Spain}\\
{\small Email: dcg@icmat.es} & {\small Email: fgancedo@us.es}
\end{tabular}
\end{quote}
\end{document}